\documentclass[12pt,a4paper]{amsart}

\usepackage{euscript,amsfonts,amssymb,amsmath,amscd,amsthm,enumerate,hyperref}

\usepackage{bbm}

\usepackage{comment}
\usepackage{todonotes}

\usepackage{mathrsfs}

\usepackage{graphicx}

\usepackage{color}

\usepackage[margin=1.1in]{geometry}

\usepackage{bbm}

\newtheorem{theorem}{Theorem}[section]
\newtheorem*{theorem*}{Theorem}
\newtheorem{lemma}[theorem]{Lemma}
\newtheorem{definition}[theorem]{Definition}
\newtheorem{proposition}[theorem]{Proposition}

\newtheorem{corollary}[theorem]{Corollary}
\newtheorem{assumption}[theorem]{Assumption}

\theoremstyle{remark}
\newtheorem{rmk}[theorem]{Remark}

\newcommand{\eps}{\varepsilon}

\newcommand{\E}{\mathbb E}

\newcommand{\A}{\mathcal A}

\newcommand{\pp}{\mathbb{P}}

\newcommand{\rr}{\mathbb{R}}

\newcommand{\nn}{\mathbb{N}}
\newcommand{\ep}{\hfill \ensuremath{\Box}}
\newcommand{\eq}{\begin{equation}}
\newcommand{\en}{\end{equation}}

\numberwithin{equation}{section}

\title[Dynamics of observables]{Dynamics of observables in rank-based models and performance of functionally generated portfolios} 

\author{Sergio A. Almada Monter}
\email{saalm56@gmail.com}

\author{Mykhaylo Shkolnikov}
\address{ORFE Department, Princeton University, Princeton, NJ 08544, USA}
\email{mshkolni@gmail.com}
\thanks{M. Shkolnikov was partially supported by the NSF grant DMS-1506290.}

\author{Jiacheng Zhang}
\address{ORFE Department, Princeton University, Princeton, NJ 08544, USA}
\email{jiacheng@princeton.edu}

\begin{document}

\begin{abstract}
In the seminal work \cite{Fe}, several macroscopic market observables have been introduced, in an attempt to find characteristics capturing the diversity of a financial market. Despite the crucial importance of such observables for investment decisions, a concise mathematical description of their dynamics has been missing. We fill this gap in the setting of rank-based models and expect our ideas to extend to other models of large financial markets as well. The results are then used to study the performance of multiplicatively and additively functionally generated portfolios, in particular, over short-term and medium-term horizons. 
\end{abstract}



\maketitle

\section{Introduction}\label{sec:intro}

A key characteristic of an equity market is its diversity which, on an intuitive level, describes how evenly the investors distribute their capital among the publicly traded companies. In the seminal work \cite{Fe} (see also \cite{Fe2}, \cite{FK}, \cite{FKK}), \textsc{Fernholz} has initiated the program of capturing the concept of diversity mathematically and, thus, quantifying its implications on the performance of investment portfolios. Given the market weights $\mu_1(t),\,\mu_2(t),\,\ldots,\,\mu_n(t)$ at a time $t\ge0$ (i.e. the fractions of market capital invested in the $n$ publicly traded companies at that time), he suggested to measure the market diversity by 
\begin{equation}\label{eq:diversity_def}
D_p(t):=\bigg(\sum_{i=1}^n \mu_i(t)^p\bigg)^{1/p}\quad\text{for some}\quad p\in(0,1)
\end{equation}
(see \cite[Example 3.4.4]{Fe}). The choice $p\in(0,1)$ ensures that the right-hand side of \eqref{eq:diversity_def} is a concave function of $\mu_1(t),\,\mu_2(t),\,\ldots,\,\mu_n(t)$ and attains its maximum for the uniform capital distribution $\mu_1(t)=\mu_2(t)=\cdots=\mu_n(t)=\frac{1}{n}$. The limiting case 
\begin{equation}\label{eq:entropy_def}
H(t):=\lim_{p\uparrow 1} \frac{p}{1-p}\log D_p(t)=-\sum_{i=1}^n \mu_i(t)\log\mu_i(t),
\end{equation}  
known as the market entropy, retains the latter two properties and can therefore be regarded as an alternative measure of the market diversity (cf. \cite[Section 2.3]{Fe}). We refer to \cite[Figures 6.7, 7.3, and 6.2]{Fe} for plots of the process $D_{1/2}(\cdot)$ for the largest 1000 companies in the U.S., the process $D_{0.76}(\cdot)$ for the companies forming the S\&P 500 index, and the process $H(\cdot)$ for the companies in the Center for Research in Securities Prices (CRSP) database of major U.S. stock exchanges (including the NYSE, the AMEX and the NASDAQ), respectively.

\medskip

Despite the considerable interest in the quantities $D_p(\cdot)$, $p\in(0,1)$ and $H(\cdot)$ (and the associated functionally generated portfolios, see below), a concise mathematical description of their dynamics has been missing so far. The challenge lies thereby in the fact that, while the vector of market weights $(\mu_1(\cdot),\mu_2(\cdot),\ldots,\mu_n(\cdot))$ is typically modeled by a Markov process, the Markov property is generally not inherited by $D_p(\cdot)$, $p\in(0,1)$ or $H(\cdot)$. Our first main goal in this paper is to capture the dynamics of non-linear macroscopic observables of the point process of logarithmic market capitalizations, such as $D_p(\cdot)$ and $H(\cdot)$, in the context of rank-based (a.k.a. first-order) models and when the number of companies $n$ is large. We choose to work with rank-based models because they are known to form the simplest class of market models that is able to reproduce the true long-term average capital distribution of a financial market (see \cite[Chapter 5]{Fe}, \cite{BFK} and \cite[Chapter 13]{FK}). We point to \cite[Figure 5.1]{Fe} for a plot of the latter for the stocks in the CRSP database. It is worth stressing that, even though the details of our proofs rely on the specifics of rank-based models, the high-level ideas of our work can be applied to any model of a large financial market.  

\medskip

The term \textit{rank-based model} refers to the unique weak solution of the system of stochastic differential equations (SDEs)
\begin{equation}\label{intro_equ1}
\mathrm{d}X_i^{(n)}(t)=b\big(F_{\varrho^{(n)}(t)}\big(X_i^{(n)}(t)\big)\big)\,\mathrm{d}t+\sigma\big(F_{\varrho^{(n)}(t)}\big(X_i^{(n)}(t)\big)\big)\,\mathrm{d}B^{(n)}_i(t),\quad i=1,\,2,\,\ldots,\,n,
\end{equation}
with coefficient functions $b:\,[0,1]\to\rr$ and $\sigma:\,[0,1]\to(0,\infty)$, the empirical cumulative distribution functions $F_{\varrho^{(n)}(t)}(x):=\frac{1}{n}\sum_{i=1}^n \mathbf{1}_{\{X^{(n)}_i(t)\le x\}}$, and independent standard Brownian motions $B^{(n)}_1,\,B^{(n)}_2,\,\ldots,\,B^{(n)}_n$. The system \eqref{intro_equ1} is a special case of the systems of SDEs studied by \textsc{Bass} and \textsc{Pardoux} in \cite{BP}, who were motivated by the piecewise linear filtering problem. In particular, the main result of \cite{BP} shows the weak uniqueness for \eqref{intro_equ1} (the weak existence for \eqref{intro_equ1} falls under the classical result of \cite[Exercise 12.4.3]{SV}). More recently, the interacting particle system described by \eqref{intro_equ1} and its variants have attracted much attention due to their appearance in stochastic portfolio theory and an open problem of \textsc{Aldous} \cite{Ald} (see \cite[Section 5.5]{Fe}, \cite{BFK}, \cite[Section 13]{FK}, \cite{IPBKF}, \cite{Sh}, \cite{JR2}, \cite{IPS}, \cite{JR}, \cite{Re} for the former and \cite{McSh}, \cite{PP}, \cite{DT}, \cite{TT}, \cite{ST}, \cite{CDSS}, \cite{DJO} for the latter). 

\medskip

The $n\to\infty$ asymptotics of non-linear macroscopic observables that we derive herein rely on the law of large numbers for rank-based models in \cite[Corollary 2.13]{JR2} (see also \cite[Corollary 1.6]{DSVZ}, \cite[Theorem 1.2]{Sh}) and the associated central limit theorem in \cite[Theorem 1.2]{KoSh}. Both of these results hold under the following (stronger than the original) assumption.

\begin{assumption}\label{main_asmp}
\begin{enumerate}[(a)]
\item There exists a probability measure $\lambda$ on $\rr$ possessing a bounded density function and satisfying 
\begin{equation}
\forall\,\theta>0:\quad \int_\rr e^{\theta|x|}\,\lambda(\mathrm{d}x)<\infty
\end{equation}
such that the initial locations of the particles $X^{(n)}_1(0),\,X^{(n)}_2(0),\,\ldots,\,X^{(n)}_n(0)$ are i.i.d. according to $\lambda$ for all $n\in\nn$.
\item The coefficient functions $b$ and $\sigma$ are differentiable with locally H\"older continuous derivatives.
\end{enumerate}
\end{assumption}

We now state the versions of \cite[Corollary 2.13]{JR2} and \cite[Theorem 1.2]{KoSh} used in this paper for future reference. Hereby, we write $M_1(\rr)$ for the space of probability measures on $\rr$ equipped with the topology of weak convergence, $C([0,\infty),M_1(\rr))$ for the space of continuous functions from $[0,\infty)$ to $M_1(\rr)$ endowed with the topology of locally uniform convergence, as well as $M_{\mathrm{fin}}(\rr)$ and $M_{\mathrm{fin}}([0,t]\times\rr)$, $t>0$ for the spaces of finite signed measures on $\rr$ and $[0,t]\times\rr$, $t>0$ viewed as the duals of $C_0(\rr)$ and $C_0([0,t]\times\rr)$, $t>0$ with the associated weak-$*$ topologies, respectively.

\begin{proposition}[cf. \cite{JR2}, Corollary 2.13] \label{prop1}
Under Assumption \ref{main_asmp} the processes of empirical measures $\varrho^{(n)}(\cdot)=\frac{1}{n}\sum_{i=1}^n \delta_{X^{(n)}_i(\cdot)}$, $n\in\nn$ converge in probability in $C([0,\infty),M_1(\rr))$ to a deterministic limit $\varrho(\cdot)$. Moreover, the corresponding process of cumulative distribution functions $R(t,\cdot):=F_{\varrho(t)}(\cdot)$, $t\ge0$ forms the unique generalized solution of the Cauchy problem for the porous medium equation
\begin{equation}\label{prop1_equ1}
R_t=-B(R)_x+\Sigma(R)_{xx},\quad R(0,\cdot)=F_{\lambda}(\cdot)
\end{equation}
in the sense of \cite[Definition 3]{Gi}, where $B(r):=\int_0^r b(a)\,\mathrm{d}a, \Sigma(r):=\int_0^r \frac{1}{2}\,\sigma(a)^2\,\mathrm{d}a$, and $F_\lambda(\cdot)$ is the cumulative distribution function of $\lambda$.
\end{proposition}

\begin{proposition}[cf. \cite{KoSh}, Theorem 1.2]\label{prop2}
Let Assumption \ref{main_asmp} be satisfied and $G$ be the mild solution of the stochastic partial differential equation (SPDE)
\begin{equation}\label{SPDE}
G_t=-\big(b(R)\,G\big)_x+\bigg(\frac{\sigma(R)^2}{2}\,G\bigg)_{xx}
+\sigma(R)\,R_x^{1/2}\,\dot{W},\quad G(0,\cdot)=\beta(F_\lambda(\cdot)),
\end{equation}
with the function $R$ from Proposition \ref{prop1}, the space-time white noise $\dot{W}$, and a standard Brownian bridge $\beta$ independent of $\dot{W}$. In other words, 
\begin{equation}\label{mild_def}
\begin{split}
G(t,x)=\int_\rr \beta(F_\lambda(y))\,p(0,y;t,x)\,\mathrm{d}y
+\int_0^t \!\int_\rr \sigma(R(s,y))\,R_x(s,y)^{1/2}\,p(s,y;t,x)\,\mathrm{d}W(s,y), \\
(t,x)\in[0,\infty)\times\rr,
\end{split}
\end{equation}
where $p$ is the transition kernel associated with the solution of the martingale problem for the operators $b(R(t,\cdot))\,\frac{\mathrm{d}}{\mathrm{d}x}+\frac{\sigma(R(t,\cdot))^2}{2}\,\frac{\mathrm{d}^2}{\mathrm{d}x^2}$, $t\ge0$ and the double integral is taken in the It\^o sense. Then, the sequences of processes
\begin{equation}\label{seq_fluct}
t\mapsto\!\sqrt{n}\big(F_{\varrho^{(n)}(t)}(x)-R(t,x)\big) \mathrm{d}x,\,n\in\nn \;\;\text{and}\;\;t\mapsto\!\sqrt{n}\big(F_{\varrho^{(n)}(s)}(x)-R(s,x)\big)\mathrm{d}x\,\mathrm{d}s,\,n\in\nn,
\end{equation}
with values in $M_{\mathrm{fin}}(\rr)$ and $M_{\mathrm{fin}}([0,t]\times\rr)$, $t\!>\!0$, respectively, converge jointly in the finite-dimensional distribution sense to 
\begin{equation}
t\mapsto G(t,x)\,\mathrm{d}x\quad\text{and}\quad t\mapsto G(s,x)\,\mathbf{1}_{[0,t]\times\rr}(s,x)\,\mathrm{d}s\,\mathrm{d}x.
\end{equation}
\end{proposition}

We are now ready to give the first two main results of the present work, which yield a comprehensive description of the large $n$ asymptotic dynamics for non-linear macroscopic observables of the form
\begin{equation}\label{intro_equ3}
{\mathcal J}_{J;f_1,\ldots,f_k}\big(\alpha(\cdot)\big):= 
J\bigg(\int_\rr f_1\,\mathrm{d}\alpha(\cdot),\,\ldots,\,\int_\rr f_k\,\mathrm{d}\alpha(\cdot)\bigg),
\end{equation}
where $\alpha(\cdot)\in C([0,\infty),M_1(\rr))$, $J$ is a continuously differentiable function, and 
\begin{equation}
f_1,\,\ldots,\,f_k\in{\mathcal E}_\ell:=\bigg\{f\in C^\ell(\rr):\;\bigg|\frac{\mathrm{d}^\ell f}{\mathrm{d}x^\ell}(x)\bigg|\le C e^{C|x|},\,x\in\rr\;\;\text{for some}\;\;C\ge0\bigg\},
\end{equation}
with $\ell=1$ in Theorem \ref{thm1} and $\ell=3$ in Theorem \ref{thm2} below. For simplicity, we use henceforth the bilinear form notation $\langle f,\nu \rangle$ for $\int_\rr f\,\mathrm{d}\nu$ and write $(f_1,\,\ldots,\,f_k)\in\mathcal{E}^U_\ell$ for a $U\subseteq\rr^k$ whenever $f_1,\,\ldots,\,f_k\in\mathcal{E}_\ell$ and
\begin{equation*}
\big(\langle f_1,\nu\rangle,\,\ldots,\,\langle f_k,\nu\rangle\big)\in U\text{ for all }\nu\in M_1(\rr)\text{ fulfilling }\int_\rr e^{\theta|x|}\,\nu(\mathrm{d}x)<\infty,\;\theta>0.    
\end{equation*} 

\begin{theorem}\label{thm1}
Let Assumption \ref{main_asmp} be satisfied and $(f_1,\,\ldots,\,f_k) \in{\mathcal E}^U_1$ for a convex open $U\subseteq\rr^k$. Then, for all $J\in C^1(U)$, one has
\begin{equation}\label{thm1_result}
\sqrt{n}\Big({\mathcal J}_{J;f_1,\ldots,f_k}\big(\varrho^{(n)}(\cdot)\big)
-\mathcal{J}_{J;f_1,\ldots,f_k}\big(\varrho(\cdot)\big)\Big)
\stackrel{n\to\infty}\longrightarrow -\sum_{j=1}^k 
{\mathcal J}_{J_{x_j};f_1,\ldots,f_k}\big(\varrho(\cdot)\big)\,\int_\rr f'_j(x)\,G(\cdot,x)\,\mathrm{d}x
\end{equation}
in the finite-dimensional distribution sense.
\end{theorem}

\begin{theorem}\label{thm2}
Let Assumption \ref{main_asmp} be satisfied, $(f_1,\,\ldots,\,f_k) \in{\mathcal E}^U_3$ for a convex open $U\subseteq\rr^k$, and $J\in C^1(U)$. Suppose $a\in\rr$ is such that 
\begin{equation}\label{asmp_on_a}
\tau:=\inf\big\{t\ge0:\,\mathcal{J}_{J;f_1,\ldots,f_k}\big(\varrho(t)\big)=a\big\}<\infty\quad\text{and}\quad \frac{\mathrm{d}\mathcal{J}_{J;f_1,\ldots,f_k}(\varrho(\cdot))}{\mathrm{d}t}(\tau)\neq 0. 
\end{equation}
Then, the sequence of hitting times
\begin{equation}
\tau^{(n)}:=\inf\Big\{t\ge0:\,{\mathcal J}_{J;f_1,\ldots,f_k}\big(\varrho^{(n)}(t)\big)=a\Big\},\quad n\in\nn
\end{equation}  
converges in distribution when properly rescaled:
\begin{equation}\label{thm2_result}
\sqrt{n}\,(\tau^{(n)}-\tau)\stackrel{n\to\infty}\longrightarrow \frac{\sum_{j=1}^k 
{\mathcal J}_{J_{x_j};f_1,\ldots,f_k}(\varrho(\tau))\,\int_\rr f'_j(x)\,G(\tau,x)\,\mathrm{d}x}{\frac{\mathrm{d}\mathcal{J}_{J;f_1,\ldots,f_k}(\varrho(\cdot))}{\mathrm{d}t}(\tau)}\,.
\end{equation}
\end{theorem}

\begin{rmk}
We emphasize that Theorems \ref{thm1} and \ref{thm2} apply, in particular, to the (appropriately normalized) processes $D_p(\cdot)$, $p\in(0,1)$ and $H(\cdot)$ of \eqref{eq:diversity_def} and \eqref{eq:entropy_def}, respectively, allowing to approximate them by Gaussian processes (Theorem \ref{thm1}) and their hitting times by Gaussian random variables (Theorem \ref{thm2}). For more details, please see Section~\ref{sec:observable} below.
\end{rmk}

\begin{rmk}
The condition of finiteness of all exponential moments on $\lambda$ in Assumption \ref{main_asmp}(a) enters naturally in the context of market observables from stochastic portfolio theory, which often (e.g. in the case of $D_p(\cdot)$, $p\in(0,1)$) involve powers of the market capitalizations that, in turn, are images of $X^{(n)}_1(\cdot),\,X^{(n)}_2(\cdot),\,\ldots,\,X^{(n)}_n(\cdot)$ under the exponential function.   
\end{rmk}

Theorem \ref{thm2} can be used further to get estimates on the performance of multiplicatively and additively functionally generated portfolios in the sense of \cite[Chapter 3]{Fe} and \cite{KaRu}, respectively. Consider a function
\begin{equation} \label{eqn:JPsi}
\Psi:\,\rr^n\to\rr,\quad x\mapsto J\bigg(\frac{1}{n}\sum_{i=1}^n f_1(x_i),\,\ldots,\,\frac{1}{n}\sum_{i=1}^n f_k(x_i)\bigg)
\end{equation}
with the homogeneity property
\begin{equation}
\forall\,x\in\rr^n,\,r\in\rr:\quad \Psi(x)=\Psi(x_1+r,x_2+r,\ldots,x_n+r).
\end{equation}
The latter is equivalent to the existence of the representation
\begin{equation}\label{eqn:PsitildePsi}
\Psi(x)=\widetilde{\Psi}\bigg(\frac{e^{x_1}}{\sum_{i=1}^n e^{x_i}},\frac{e^{x_2}}{\sum_{i=1}^n e^{x_i}},\ldots,\frac{e^{x_n}}{\sum_{i=1}^n e^{x_i}}\bigg),\quad x\in\rr^n.
\end{equation}  
If $\widetilde{\Psi}$ can be extended to a twice continuously differentiable function on an open neighborhood of the open unit simplex $\{x\in(0,1)^n:\,\sum_{i=1}^n x_i=1\}\subset\rr^n$, then one can formally define the weights (i.e. the fractions invested in the different companies) of the portfolios $\pi^{\widetilde{\Psi};\times}$ and $\pi^{\widetilde{\Psi};+}$ multiplicatively and additively generated by $\widetilde{\Psi}$ via
\begin{eqnarray*}
&&\;\;\;\pi^{\widetilde{\Psi};\times}_i(\cdot)=\Big(\!(\log \widetilde{\Psi})_{x_i}\big(\mu(\cdot)\big)+1-\sum_{j=1}^n \mu_j(\cdot)\,(\log \widetilde{\Psi})_{x_j}\big(\mu(\cdot)\big)\!\Big)\mu_i(\cdot),\;i\!=\! 1,2,\ldots,n, \\
&&\;\;\;\pi^{\widetilde{\Psi};+}_i(\cdot)=\bigg(\!\frac{\widetilde{\Psi}_{x_i}(\mu(\cdot))
\!-\!\sum_{j=1}^n \mu_j(\cdot)\,\widetilde{\Psi}_{x_j}(\mu(\cdot))}{\widetilde{\Psi}(\mu(\cdot))\!-\!\frac{1}{2}\sum_{i,j=1}^n \!\int_0^\cdot \widetilde{\Psi}_{x_ix_j}(\mu(s)) \mathrm{d}[\mu_i,\mu_j](s)}\!+\! 1\!\bigg)\mu_i(\cdot),\;i\!=\! 1,2,\ldots,n,
\end{eqnarray*}
respectively. Here, $\mu(\cdot)=(\mu_1(\cdot),\mu_2(\cdot),\ldots,\mu_n(\cdot))$ is the process of the market weights, which in the context of a rank-based model amounts to
\begin{equation*}
\bigg(\frac{e^{X_1(\cdot)}}{\sum_{i=1}^n e^{X_i(\cdot)}},\frac{e^{X_2(\cdot)}}{\sum_{i=1}^n e^{X_i(\cdot)}},\ldots,
\frac{e^{X_n(\cdot)}}{\sum_{i=1}^n e^{X_i(\cdot)}}\bigg),
\end{equation*}
and $[\cdot,\cdot]$ denotes the quadratic covariation process. 

\medskip

Functionally generated portfolios $\pi^{\widetilde{\Psi};\times}(\cdot)$, $\pi^{\widetilde{\Psi};+}(\cdot)$ have the remarkable property that their values $V^{\widetilde{\Psi};\times}(\cdot)$, $V^{\widetilde{\Psi};+}(\cdot)$ relative to that of the market portfolio $\mu(\cdot)$ admit pathwise representations, which under the usual convention $V^{\widetilde{\Psi};\times}(0)=V^{\widetilde{\Psi};+}(0)=1$ read 
\begin{eqnarray}
&& V^{\widetilde{\Psi};\times}(t)=\frac{\widetilde{\Psi}(\mu(t))}{\widetilde{\Psi}(\mu(0))}\,\exp\bigg(\!-\frac{1}{2}\sum_{i,j=1}^n \int_0^t \frac{\widetilde{\Psi}_{x_ix_j}(\mu(\cdot))}{\widetilde{\Psi}(\mu(\cdot))}\,\mathrm{d}[\mu_i,\mu_j](\cdot)\bigg),\;\; t\ge 0, \label{master_mult} \\
&& V^{\widetilde{\Psi};+}(t)=1+\widetilde{\Psi}\big(\mu(t)\big)-\widetilde{\Psi}\big(\mu(0)\big)-\frac{1}{2}\sum_{i,j=1}^n \int_0^t \widetilde{\Psi}_{x_ix_j}\big(\mu(\cdot)\big)\,\mathrm{d}[\mu_i,\mu_j](\cdot),\;\; t\ge0 \label{master_add}
\end{eqnarray}
(cf. \cite[equation (11.2)]{FK}, \cite[equation (4.3)]{KaRu}). We assume henceforth that the function $\widetilde{\Psi}$ is positive and concave in the setting of \eqref{master_mult} or concave in the setting of \eqref{master_add}, since then the respective \textit{excess growth process} $-\frac{1}{2}\sum_{i,j=1}^n \int_0^t \frac{\widetilde{\Psi}_{x_ix_j}(\mu(\cdot))}{\widetilde{\Psi}(\mu(\cdot))}\,\mathrm{d}[\mu_i,\mu_j](\cdot)$, $t\ge0$ or $-\frac{1}{2}\sum_{i,j=1}^n \int_0^t \widetilde{\Psi}_{x_ix_j}\big(\mu(\cdot)\big)\,\mathrm{d}[\mu_i,\mu_j](\cdot)$, $t\ge0$ is non-decreasing (see \cite[Example 3.5]{KaRu}). In particular, the associated value process $V^{\widetilde{\Psi};\times}(\cdot)$ or $V^{\widetilde{\Psi};+}(\cdot)$ reaches levels $v>1$ before the hitting times $\tau^{(n)}$ of levels $a=v\,\widetilde{\Psi}(\mu(0))$ or $a=v-1+\widetilde{\Psi}(\mu(0))$, respectively, whose asymptotics are described by Theorem \ref{thm2}.

\medskip

More precise estimates on the processes $V^{\widetilde{\Psi};\times}(\cdot)$ and $V^{\widetilde{\Psi};+}(\cdot)$ can be obtained under the additional assumptions
\begin{equation}\label{extra_asmp}
\forall\,n\in\nn,\,i=1,2,\ldots,n-1:\;\frac{1}{i}\,\sum_{j=1}^i b\Big(\frac{j}{n}\Big) > \frac{1}{n-i}\,\sum_{j=i+1}^n b\Big(\frac{j}{n}\Big)\quad\text{and}\quad\sigma(\cdot)=1
\end{equation}
(see \cite[Remark on p.~2187]{PP} for a detailed discussion of the first assumption; in addition, note that the constant $1$ in the second assumption can be turned into any other positive constant by a deterministic time change). Indeed, under the assumptions in \eqref{extra_asmp}, \cite[Corollary 8]{IPS} applies and can be naturally combined with Theorem \ref{thm2}.

\begin{corollary}\label{main_cor}
Let Assumption \ref{main_asmp}, the assumptions in \eqref{extra_asmp}, and for all $n\in\nn$,
\begin{equation}\label{mult_range}
\mathrm{ess\,sup}_{\omega,t} \sum_{i,j=1}^n \frac{\widetilde{\Psi}_{x_ix_j}(\mu(\cdot))}{\widetilde{\Psi}(\mu(\cdot))}\,\frac{\mathrm{d}[\mu_i,\mu_j](\cdot)}{\mathrm{d}t} 
-\mathrm{ess\,inf}_{\omega,t} \sum_{i,j=1}^n \frac{\widetilde{\Psi}_{x_ix_j}(\mu(\cdot))}{\widetilde{\Psi}(\mu(\cdot))}\,\frac{\mathrm{d}[\mu_i,\mu_j](\cdot)}{\mathrm{d}t}\in(0,\infty)
\end{equation}
or
\begin{equation}\label{add_range}
\mathrm{ess\,sup}_{\omega,t} \sum_{i,j=1}^n \widetilde{\Psi}_{x_ix_j}(\mu(\cdot))\,\frac{\mathrm{d}[\mu_i,\mu_j](\cdot)}{\mathrm{d}t} 
-\mathrm{ess\,inf}_{\omega,t} \sum_{i,j=1}^n \widetilde{\Psi}_{x_ix_j}(\mu(\cdot))\,\frac{\mathrm{d}[\mu_i,\mu_j](\cdot)}{\mathrm{d}t}\in(0,\infty)
\end{equation}
be satisfied. Then, for the functions $J$ and $f_1,\,\ldots,\,f_k$ of \eqref{eqn:JPsi}, \eqref{eqn:PsitildePsi}, $a,\,\tau$ fulfilling the conditions in \eqref{asmp_on_a}, $r^\times:=-\underset{t\to\infty}{\lim}\,\frac{1}{2t}\int_0^t \frac{\widetilde{\Psi}_{x_ix_j}(\mu(\cdot))}{\widetilde{\Psi}(\mu(\cdot))}\,\mathrm{d}[\mu_i,\mu_j](\cdot)$ or $r^+:=-\underset{t\to\infty}{\lim}\,\frac{1}{2t}\sum_{i,j=1}^n \int_0^t \,\widetilde{\Psi}_{x_ix_j}(\mu(\cdot))\,\mathrm{d}[\mu_i,\mu_j](\cdot)$, and $r,s>0$, the stopping times
\begin{eqnarray*}
&& \eta^{\widetilde{\Psi};\times}:=\inf\Big\{t\ge0:\,V^{\widetilde{\Psi};\times}(t)=\frac{a}{\widetilde{\Psi}(\mu(0))}\,e^{(r^\times-r)(\tau-s/\sqrt{n})}\Big\}, \\
&& \eta^{\widetilde{\Psi};+}:=\inf\Big\{t\ge0:\,V^{\widetilde{\Psi};+}(t)=1+a-\widetilde{\Psi}\big(\mu(0)\big)+(r^+-r)\big(\tau-s/\sqrt{n}\big)\Big\}
\end{eqnarray*}
satisfy for all $\eps>0$ the respective estimates
\begin{eqnarray}
&& \;\,\pp\big(\eta^{\widetilde{\Psi};\times}\ge \tau+s/\sqrt{n}\big)\le 2\overline{\Phi}(s/\chi)\big(1+o_n(1)\big)+\bigg\|\frac{\mathrm{d}\kappa^{(n)}}{\mathrm{d}\zeta^{(n)}}
\bigg\|_{L^2(\zeta^{(n)})}
e^{-c^\times(r,\eps)(\tau-s/\sqrt{n})}, \label{mult_port_bnd} \\
&& \;\,\pp\big(\eta^{\widetilde{\Psi};+}\ge \tau+s/\sqrt{n}\big)\le 2\overline{\Phi}(s/\chi)\big(1+o_n(1)\big)+\bigg\|\frac{\mathrm{d}\kappa^{(n)}}{\mathrm{d}\zeta^{(n)}}\bigg\|_{L^2(\zeta^{(n)})}e^{-c^+(r,\eps)(\tau-s/\sqrt{n})}. \label{add_port_bnd}
\end{eqnarray}
Hereby, $\overline{\Phi}$ is the standard normal tail cumulative distribution function; $\chi$ is the standard deviation of the random variable on the right-hand side of \eqref{thm2_result}; $o_n(1)$ is a quantity tending to $0$ as $n\to\infty$; $\kappa^{(n)}$ and $\zeta^{(n)}$ are the laws of the vector of differences between the consecutive order statistics of $(X^{(n)}_1,X^{(n)}_2,\ldots,X^{(n)}_n)$ at time $0$ and in stationarity, respectively; and
\begin{equation} \label{eqn: c}
\begin{split}
c^\times(r,\eps)=&\,\frac{\underset{1\le j\le n-1}{\min}\big(\sum_{i=1}^j b(\frac{i}{n})-\frac{j}{n}\sum_{i=1}^n b(\frac{i}{n})\big)^2}{2-2\cos\frac{\pi}{n}} \\
&\,\cdot\max\Bigg(\frac{r^2}{(C^\times)^2},4\eps(\eps+v^\times)\Bigg(\sqrt{1+\frac{r^2}{2\eps(\eps+v^\times)^2\max(|C^{\times,\uparrow}|,|C^{\times,\downarrow}|)^2}}-1\Bigg)\!\Bigg), 
\end{split}
\end{equation}
\begin{equation}
\begin{split}
c^+(r,\eps)=&\,\frac{\underset{1\le j\le n-1}{\min}\big(\sum_{i=1}^j b(\frac{i}{n})-\frac{j}{n}\sum_{i=1}^n b(\frac{i}{n})\big)^2}{2-2\cos\frac{\pi}{n}} \\
&\,\cdot\max\Bigg(\frac{r^2}{(C^+)^2},4\eps(\eps+v^+)\Bigg(\sqrt{1+\frac{r^2}{2\eps(\eps+v^+)^2\max(|C^{+,\uparrow}|,|C^{+,\downarrow}|)^2}}-1\Bigg)\!\Bigg) 
\end{split}
\end{equation}
with $C^\times$, $C^{\times,\uparrow}$, $C^{\times,\downarrow}$, and $v^\times$ ($C^+$, $C^{+,\uparrow}$, $C^{+,\downarrow}$, and $v^+$ resp.) being the overall expression, the essential supremum, the essential infimum, and the variance under $\zeta^{(n)}$ of the expression inside the essential supremum in \eqref{mult_range} (\eqref{add_range} resp.).
\end{corollary}

\begin{rmk}
The inequalities \eqref{mult_port_bnd}, \eqref{add_port_bnd} can be interpreted as follows. If one invests in the portfolio generated multiplicatively (or additively resp.) by a function $\widetilde{\Psi}$ satisfying the condition \eqref{mult_range} (or \eqref{add_range} resp.) and aims for the associated process $\widetilde{\Psi}(\mu(\cdot))$ to reach an admissible value of $a$ (i.e. one for which \eqref{asmp_on_a} holds), then one will achieve a logarithmic (or arithmetic resp.) return relative to the market portfolio $\mu(\cdot)$ of $\log a-\log\widetilde{\Psi}(\mu(0))+(r^\times-r)(\tau-s/\sqrt{n})$ (or $a-\widetilde{\Psi}(\mu(0))+(r^+-r)(\tau-s/\sqrt{n})$ resp.) before time $\tau+s/\sqrt{n}$ with a confidence probability of at least one minus the right-hand side of \eqref{mult_port_bnd} (or \eqref{add_port_bnd} resp.). We note that, in the practically relevant regime $a\ge{\mathcal J}_{J;f_1,\ldots,f_k}(\varrho(0))$ and $r^\times\ge r$ (or $r^+\ge r$ resp.), by increasing the values of $a$ and $s$ one can increase both the relative return and the confidence probability for sufficiently large $n$, at the expense of thereby increasing the upper bound $\tau+s/\sqrt{n}$ on the investment horizon. 
\end{rmk}

The rest of the paper is structured as follows. In Section \ref{sec:prelim}, we collect some results from \cite{Ha}, \cite{KoSh} and \cite{BoLe} that are used repeatedly in the proofs of Theorems \ref{thm1} and \ref{thm2}. Section~\ref{sec:PfLemma3} is then devoted to the proof of Theorem \ref{thm1}. The latter is based on Proposition \ref{prop2}, but requires significant additional work due to the exponential growth at infinity of the derivatives of functions in ${\mathcal E}_1$ and the non-linearity of $J$. In particular, the proof invokes the mean stochastic comparison of \cite{Ha} and the quantitative propagation of chaos result of \cite[Theorem 1.6]{KoSh}. Subsequently, we give the proof of Theorem \ref{thm2} in Section \ref{sec:proofCLT}, which relies on the previously mentioned tools and a creative reduction to the estimate on the expected Wasserstein distance $W_1$ between the empirical measure of an i.i.d. sample and the underlying distribution in \cite[Theorem 3.2]{BoLe}. Next, in Section \ref{sec:observable}, we apply Theorems \ref{thm1} and \ref{thm2} to the main examples of diversity measures from stochastic portfolio theory. Lastly, in Section \ref{sec:functionallyGenerated}, we provide the proof of Corollary \ref{main_cor}. 

\medskip

\noindent\textbf{Acknowledgement.} The second author (M.~S.) thanks Adrian Banner, Robert Fernholz, Ioannis Karatzas and Vassilios Papathanakos for bringing the problem of describing the dynamics of the market diversity to his attention. The second and the third authors (M.~S. and J.~Z.) are grateful to the participants of the research meetings  at \textsc{Intech} Investment Management for their insightful comments.

\section{Preliminaries} \label{sec:prelim}

The general \textit{propagation of chaos} paradigm (see \cite{SzChaos}) suggests that under Assumption \ref{main_asmp}, for large values of $n$, the weak solution of \eqref{intro_equ1} should be well-approximated by the strong solution of 
\begin{equation}\label{pre_1}
d\overline{X}_i^{(n)}(t)=b\big(R(t,\overline{X}_i^{(n)}(t))\big)\,\mathrm{d}t+\sigma\big(R(t,\overline{X}_i^{(n)}(t))\big)\,\mathrm{d}B_i^{(n)}(t),\quad i=1,\,2,\,\ldots,\,n
\end{equation}
with the initial condition $\overline{X}_i^{(n)}(0)=X_i^{(n)}(0)$, $i=1,2,\ldots,n$, the function $R$ from Proposition \ref{prop1} and the same Brownian motions $B_1^{(n)},\,B_2^{(n)},\,\ldots,\,B_n^{(n)}$ as in \eqref{intro_equ1}. Indeed, under Assumption \ref{main_asmp}, the coefficient functions $(t,x)\mapsto b(R(t,x))$, $(t,x)\mapsto\sigma(R(t,x))$ are uniformly Lipschitz in $x$ on any strip of the form $[0,T]\times\rr$ by \cite[Proposition 2.5]{KoSh}, so that the strong existence and uniqueness for \eqref{pre_1} readily follow. To quantify the term ``well-approximated'' we introduce the process of empirical measures 
\begin{equation} \label{eqn: barRho}
\overline{\varrho}^{(n)}(\cdot)= \frac{1}{n}\sum_{i=1}^n\delta_{\overline{X}_i^{(n)}(\cdot)}
\end{equation}
and recall that the Wasserstein distance of order $p\geq1$ is defined for any $\nu_1,\nu_2\in M_1(\rr)$ with finite $p$-th moments by
\begin{equation} \label{eqn:wassesteinDef}
W_p(\nu_1,\nu_2)=\inf_{Z_1\stackrel{d}{=}\nu_1,Z_2\stackrel{d}{=}\nu_2}
\E\big[\vert Z_1-Z_2 \vert^p \big]^{1/p}.
\end{equation}
Then, under Assumption \ref{main_asmp}, the following quantitative propagation of chaos estimates from \cite[Theorem 1.6]{KoSh} apply. 

\begin{proposition}[cf. \cite{KoSh}, Theorem 1.6] \label{prop2.2}
Let Assumption \ref{main_asmp} be satisfied. Then, for all $p,T>0$, one can find a constant $C=C(p,T)<\infty$ such that 
\begin{equation}\label{prop_of_chaos1}
\E\Big[\sup_{t\in[0,T]}|X_i^{(n)}(t)-\overline{X}_i^{(n)}(t)|^p\Big]\leq Cn^{-p/2}, \quad i=1,\,2,\,\ldots,\,n,\;\;n\in\nn.
\end{equation}
In particular, for $p\ge1$, it holds
\begin{equation}
\E\Big[\sup_{t\in[0,T]}W_p\big(\varrho^{(n)}(t),\overline{\varrho}^{(n)}(t)\big)\Big]\leq Cn^{-p/2},\quad n\in\nn.
\end{equation}
\end{proposition}

\smallskip

Moreover, each $\overline{\varrho}^{(n)}(t)$ constitutes the empirical measure of an i.i.d. sample from the probability measure $\varrho(t)$ introduced in Proposition \ref{prop1}. Hence, we may aim to bound the associated expected $W_1$-distance $\E[W_1(\overline{\varrho}^{(n)}(t),\varrho(t))]$ by means of \cite[Theorem 3.2]{BoLe}, which requires a moment estimate for $\varrho(t)$. The latter, in turn, can be obtained under Assumption \ref{main_asmp} from the mean stochastic comparison results of \cite{Ha} as follows. With $C^\uparrow_b:=\max_{a\in[0,1]} b(a)$, $C^\downarrow_b:=\min_{a\in[0,1]} b(a)$ and $C_\sigma:=\max_{a\in[0,1]} |\sigma(a)|$, consider the Brownian motions
\begin{eqnarray}
&& \mathrm{d}Y^\uparrow(t)=C^\uparrow_b\,\mathrm{d}t+C_\sigma\,\mathrm{d}B^{(1)}_1(t),\quad Y^\uparrow(0)\stackrel{d}{=}\lambda, \\
&& \mathrm{d}Y^\downarrow(t)=C^\downarrow_b\,\mathrm{d}t+C_\sigma\,\mathrm{d}B^{(1)}_1(t),\quad Y^\downarrow(0)\stackrel{d}{=}\lambda.
\end{eqnarray}
The next proposition is then a direct consequence of \cite[inequality (1.5) and p.~318, Remark (4)]{Ha}.

\begin{proposition}\label{prop2.4}
Let Assumption \ref{main_asmp} be satisfied. Then, for all $i=1,2,\ldots,n$, $n\in\nn$, $T>0$, $M\in\rr$ and $\theta>0$, one has the comparison results
\begin{equation}\label{Pcomp}
\begin{split}
\pp\Big(\sup_{t\in[0,T]} \big|X_i^{(n)}(t)\big|\geq M\Big)\vee 
\pp\Big(\sup_{t\in[0,T]} \big|\overline{X}_i^{(n)}(t)\big|\geq M\Big) 
\leq &\,\,2 \pp\Big(\sup_{t\in[0,T]} Y^\uparrow(t)\geq M\Big) \\
& + 2\pp\Big(\sup_{t\in[0,T]} (-Y^\downarrow(t))\geq M\Big),
\end{split}
\end{equation}
\begin{equation}\label{Ecomp}
\!\E\Big[\sup_{t\in[0,T]} e^{\theta|X_i^{(n)}(t)|}\Big]\!\vee
\E\Big[\sup_{t\in[0,T]} e^{\theta|\overline{X}_i^{(n)}(t)|}\Big] 
\!\leq\!\E\Big[\sup_{t\in[0,T]} e^{\theta Y^\uparrow(t)}\Big]
\!+\E\Big[\sup_{t\in[0,T]} e^{-\theta Y^\downarrow(t)}\Big]\!\!<\!\infty. 
\end{equation}
\end{proposition}

\smallskip

In particular, the bound of \eqref{Ecomp} allows us to use \cite[Theorem 3.2]{BoLe} to estimate each of the quantities $\E[W_1(\overline{\varrho}^{(n)}(t),\varrho(t))]$. Hereby, we keep in mind the alternative representation of the $W_1$-distance as the $L^1$-distance between the cumulative distribution functions (see e.g. \cite[Theorem 2.9]{BoLe}):
\begin{equation}\label{W1repr}
W_1(\nu_1,\nu_2)=\int_\rr \big|F_{\nu_1}(x)-F_{\nu_2}(x)\big|\,\mathrm{d}x. 
\end{equation}

\begin{proposition}\label{prop_2.3}
Let Assumption \ref{main_asmp} be satisfied. Then, for all $T>0$, one can find a constant $C=C(T)<\infty$ such that
\begin{equation}
\sup_{t\in[0,T]} \E\big[W_1\big(\overline{\varrho}^{(n)}(t),\varrho(t)\big)\big]\le Cn^{-1/2},\quad n\in\nn.
\end{equation}
\end{proposition}

\section{Proof of Theorem \ref{thm1}} \label{sec:PfLemma3}

Our starting point for the proof of Theorem \ref{thm1} is the identity
\begin{equation}\label{thm1pf_equ1}
\begin{split}
& \;\sqrt{n}\,\Big({\mathcal J}_{J;f_1,\ldots,f_k}\big(\varrho^{(n)}(\cdot)\big) - {\mathcal J}_{J;f_1,\ldots,f_k}\big(\varrho(\cdot)\big)\Big) \\
&= \sqrt{n}\,\Big(\big\langle f_1,\varrho^{(n)}(\cdot)-\varrho(\cdot)\big\rangle,\ldots,\big\langle f_k,\varrho^{(n)}(\cdot)-\varrho(\cdot)\big\rangle\Big) 
\,\nabla J\Big(\big\langle f_1,\widetilde{\varrho}^{(n)}(\cdot)\big\rangle,\ldots,\big\langle f_k,\widetilde{\varrho}^{(n)}(\cdot)\big\rangle\Big)
\end{split}
\end{equation}
due to the mean value theorem, where $\widetilde{\varrho}^{(n)}(\cdot)=\xi^{(n)}(\cdot)\varrho^{(n)}(\cdot)+(1-\xi^{(n)}(\cdot))\varrho(\cdot)$ and $\xi^{(n)}(\cdot)$ can be chosen as stochastic processes with values in $M_1(\rr)$ and $[0,1]$, respectively, by the Borel selection result of \cite[Theorem 6.9.6]{Bog}. The proof of Theorem \ref{thm1} is carried out by studying the convergence of the vector-valued stochastic processes  
\begin{eqnarray}
&& I_1^{(n)}(\cdot):=\sqrt{n}\Big(\big\langle f_1,\varrho^{(n)}(\cdot)-\varrho(\cdot)\big\rangle,\,\ldots,\,\big\langle f_k,\varrho^{(n)}(\cdot)-\varrho(\cdot)\big\rangle\Big), \label{I1def} \\
&& I_2^{(n)}(\cdot):=\Big(\big\langle f_1,\widetilde{\varrho}^{(n)}(\cdot) \big\rangle,\,\ldots,\,\big\langle f_k,\widetilde{\varrho}^{(n)}(\cdot)\big\rangle\Big) \label{I2def}
\end{eqnarray}
as $n\to\infty$. In both cases, it is helpful to introduce, for each $M>0$, an auxiliary function $h_M\in C^\infty(\rr)$ with values in $[0,1]$ such that $h_M(x)=1$ if $|x|\le M$, $h_M(x)=0$ if $|x|>M+1$, and 
\begin{equation}
\sup_{M>0}\,\sup_{x\in\rr}\,|h_M'(x)|\,\vee\,\sup_{M>0}\,\sup_{x\in\rr}\,|h_M''(x)|<\infty.
\end{equation}
In addition, we denote $(1-h_M)$ by $\widehat{h}_M$ for each $M>0$.  

\medskip

\noindent\textbf{Convergence of $I_1^{(n)}(\cdot)$.} With the mild solution $G$ of the SPDE \eqref{SPDE}, we claim that 
\begin{equation}\label{thm1pf_equ15}
I_1^{(n)}(\cdot)\stackrel{n\to\infty}{\longrightarrow} 
\bigg(\int_\rr f'_1(x)\,G(\cdot,x)\,\mathrm{d}x,\,\ldots,\,\int_\rr f'_k(x)\,G(\cdot,x)\,\mathrm{d}x\bigg)
\end{equation} 
in the finite-dimensional distribution sense. To this end, we write each component $\sqrt{n}\,\langle f_j,\varrho^{(n)}(\cdot)-\varrho(\cdot)\rangle$ of $I_1^{(n)}(\cdot)$ as
\begin{equation}\label{4summand_dec}
\begin{split}
& \sqrt{n}\!\int_\rr\! f_j'(x)\,h_M(x) \big(F_{\rho^{(n)}(\cdot)}(x)\!-\!R(\cdot,x)\big) \mathrm{d}x 
\!+\!\sqrt{n}\!\int_\rr\! f_j(x)\,h_M'(x) \big(F_{\rho^{(n)}(\cdot)}(x)\!-\!R(\cdot,x)\big) \mathrm{d}x \\
&\!+\!\sqrt{n}\!\int_\rr\! f_j'(x)\,\widehat{h}_M(x) \big(F_{\rho^{(n)}(\cdot)}(x)\!-\! R(\cdot,x)\big) \mathrm{d}x
\!+\!\sqrt{n}\!\int_\rr\! f_j(x)\,\widehat{h}_M'(x) \big(F_{\rho^{(n)}(\cdot)}(x)\!-\!R(\cdot,x)\big) \mathrm{d}x
\end{split}
\end{equation}
using $f_j=f_j\,h_M+f_j\,\widehat{h}_M$ and integration by parts (observe that the boundary terms thereby vanish thanks to $f_j\,\widehat{h}_M\in\mathcal{E}_0$, the estimate \eqref{Ecomp} and Markov's inequality). 

\medskip

Taking first the $n\to\infty$ limit and then the $M\to\infty$ limit of the first summand in \eqref{4summand_dec} for $j=1,2,\ldots,k$ gives 
\begin{equation}\label{Gaussian_lim}
\bigg(\int_\rr f_1'(x)\,G(\cdot,x)\,\mathrm{d}x,\,\ldots,\,\int_\rr f'_k(x)\,G(\cdot,x)\,\mathrm{d}x\bigg) 
\end{equation}
in the finite-dimensional distribution sense. Indeed, by Proposition \ref{prop2} the $n\to\infty$ limit results in the mean zero Gaussian process  
\begin{equation}
\bigg(\int_\rr f_1'(x)\,h_M(x)\,G(\cdot,x)\,\mathrm{d}x,\,\ldots,\,
\int_\rr f_k'(x)\,h_M(x)\,G(\cdot,x)\,\mathrm{d}x\bigg).
\end{equation}
For its convergence in finite-dimensional distribution as $M\to\infty$ to the mean zero Gaussian process in \eqref{Gaussian_lim}, it suffices to verify the convergence of the corresponding covariance functions. Upon a decomposition of $f'_1,f'_2,\ldots,f'_k$ into the positive and negative parts, the positivity of the covariance function of $G$ (see \cite[Remark 1.4]{KoSh}) and the monotone convergence theorem allow to reduce the convergence of the covariance functions to a statement about the uniform boundedness of the variances involved. 

\begin{lemma}\label{lemma:bounded_var} 
Let Assumption \ref{main_asmp} be satisfied. Then, for all $t\ge0$ and $j\in\{1,2,\ldots,k\}$,
\begin{equation}\label{bounded_var_eq}
\sup_{M>0}\;\E\Bigg[\!\bigg(\int_\rr f_j'(x)_+\,h_M(x)\,G(t,x)\,\mathrm{d}x\bigg)^{\!2}\Bigg]\!\vee\E\Bigg[\!\bigg(\int_\rr f_j'(x)_-\,h_M(x)\,G(t,x)\,\mathrm{d}x\bigg)^{\!2}\Bigg]<\infty.
\end{equation}
\end{lemma}

Assuming Lemma \ref{lemma:bounded_var}, the proof of \eqref{thm1pf_equ15} hinges on the next lemma, which shows that the contributions of the second, third and fourth summands in \eqref{4summand_dec} to the $n\to\infty$ limit of $I_1^{(n)}(\cdot)$ become negligible as $M$ tends to infinity. 

\begin{lemma}\label{lemma:unifromIntegrability1}
Let Assumption \ref{main_asmp} be satisfied. Then, for any $\eps>0$, $t\ge0$, $f_0\in\mathcal{E}_0$ and uniformly bounded family of functions $g_M\!:\rr\to\rr$, $M>0$ such that $g_M(x)=0$, $x\in[-M,M]$ for each $M>0$, 
\begin{equation}\label{unifromIntegrability1_eq}
\limsup_{M\to\infty}\,\limsup_{n\to\infty}\;\pp\bigg(\bigg|\sqrt{n}\int_\rr f_0(x)\,g_M(x)\,\big(F_{\varrho^{(n)}(t)}(x)-R(t,x)\big)\,\mathrm{d}x\bigg|>\eps\bigg)=0. 
\end{equation}
\end{lemma}

\smallskip

We proceed to the proofs of the two lemmas. 

\medskip

\noindent\textbf{Proof of Lemma \ref{lemma:bounded_var}.} For all $M>0$, we have
\begin{equation}
\begin{split}
& \;\E\Bigg[\!\bigg(\int_\rr f_j'(x)_+\,h_M(x)\,G(t,x)\,\mathrm{d}x\bigg)^{\!2}\Bigg] \\
& \le\liminf_{n\to\infty}\;\E\Bigg[\!\bigg(\int_\rr f_j'(x)_+\,h_M(x)\,\sqrt{n}\big(F_{\varrho^{(n)}(t)}(x)-R(t,x)\big)\,\mathrm{d}x\bigg)^{\!2}\Bigg] 
\end{split} 
\end{equation}
by Proposition \ref{prop2}, Skorokhod's representation theorem and Fatou's lemma. With $f_{j,M;+}(x):=\int_0^x f_j'(y)_+\,h_M(y)\,\mathrm{d}y$, integration by parts yields for the term inside the latter limit inferior
\begin{equation}
\begin{split}
& \;n\,\E\Big[\big\langle f_{j,M;+},\varrho^{(n)}(t)-\varrho(t)\big\rangle^2\Big] \\
& \le 2n\,\E\Big[\big\langle f_{j,M;+},\varrho^{(n)}(t)-\overline{\varrho}^{(n)}(t)\big\rangle^2\Big]
+2n\,\E\Big[\big\langle f_{j,M;+},\overline{\varrho}^{(n)}(t)-\varrho(t)\big\rangle^2\Big].
\end{split}
\end{equation}
Next, we insert the definitions of $\varrho^{(n)}(t)$, $\overline{\varrho}^{(n)}(t)$, apply the Cauchy-Schwarz inequality, and exploit the independence of $\overline{X}^{(n)}_1(t)\stackrel{d}{=}\overline{X}^{(n)}_2(t)\stackrel{d}{=}\cdots\stackrel{d}{=}\overline{X}^{(n)}_n(t)\stackrel{d}{=}\varrho(t)$ to get 
\begin{equation}\label{Lemma3.1_2terms}
2\,\E\bigg[\sum_{i=1}^n \Big(f_{j,M;+}\big(X^{(n)}_i(t)\!\big)\!-\!
f_{j,M;+}\big(\overline{X}^{(n)}_i(t)\!\big)\!\Big)^{\!2}\bigg]
\!+\!2\,\E\bigg[\Big(f_{j,M;+}\big(\overline{X}^{(n)}_1(t)\!\big)\!-\!
\big\langle f_{j,M;+},\varrho(t)\big\rangle\!\Big)^{\!2}\bigg].
\end{equation}
Since $f_1,\,\ldots,\,f_k\in\mathcal{E}_1$, we can pick a constant $C<\infty$ independent of $j$ and $M$ such that $|f'_{j,M;+}(x)|\le Ce^{C|x|}$, $x\in\rr$ and $|f_{j,M;+}(x)|\le Ce^{C|x|}$, $x\in\rr$. This, the convexity of the absolute value function, and the observation $(X^{(n)}_1(t),\overline{X}^{(n)}_1(t))\stackrel{d}{=}(X^{(n)}_2(t),\overline{X}^{(n)}_2(t))\stackrel{d}{=}\cdots\stackrel{d}{=}(X^{(n)}_n(t),\overline{X}^{(n)}_n(t))$ allow to bound the expression in \eqref{Lemma3.1_2terms} from above by
\begin{equation}\label{Lemma3.1_3terms}
\begin{split}
&\; 2nC^2\,\E\Big[\mathbf{1}_{\{|X^{(n)}_1(t)|\ge|\overline{X}^{(n)}_1(t)|\}} \,e^{2C|X^{(n)}_1(t)|} \big(X^{(n)}_1(t)-\overline{X}^{(n)}_1(t)\big)^{\!2}\Big] \\
& +2nC^2\,\E\Big[\mathbf{1}_{\{|\overline{X}^{(n)}_1(t)|>|X^{(n)}_1(t)|\}}\,e^{2C|\overline{X}^{(n)}_1(t)|} \big(X^{(n)}_1(t)-\overline{X}^{(n)}_1(t)\big)^{\!2}\Big]
+2C^2\,\E\Big[e^{2C|\overline{X}^{(n)}_1(t)|}\Big].
\end{split}
\end{equation}
By dropping the indicator random variables, using the Cauchy-Schwarz inequality twice, and invoking the estimate \eqref{Ecomp} and the $p=4$ version of the inequality \eqref{prop_of_chaos1} we conclude that the quantity in \eqref{Lemma3.1_3terms} is uniformly bounded in $n$ and $M$. An analogous argument for the second expectation in \eqref{bounded_var_eq} completes the proof of the lemma. \ep

\medskip

\noindent\textbf{Proof of Lemma \ref{lemma:unifromIntegrability1}.} With $F_{0,M}(x):=\int_0^x f_0(y)\,g_M(y)\,\mathrm{d}y$, we integrate by parts to rewrite the probability in \eqref{unifromIntegrability1_eq} as
\begin{equation}\label{2probab_est}
\begin{split}
& \;\pp\Big(\big|\sqrt{n}\,\big\langle F_{0,M},\varrho^{(n)}(t)-\varrho(t)\big\rangle\big|>\eps\Big) \\
& \le \pp\Big(\big|\sqrt{n}\,\big\langle F_{0,M},\varrho^{(n)}(t)-\overline{\varrho}^{(n)}(t)\big\rangle\big|>\eps/2\Big)
+ \pp\Big(\big|\sqrt{n}\,\big\langle F_{0,M},\overline{\varrho}^{(n)}(t)-\varrho(t)\big\rangle\big|>\eps/2\Big)
\end{split}
\end{equation} 
(note that the boundary terms in the integration by parts vanish thanks to $F_{0,M}\in\mathcal{E}_0$ and the estimate \eqref{Ecomp} in conjunction with Markov's inequality). 

\medskip

Now, we employ Markov's inequality, plug in the definitions of $\varrho^{(n)}(t)$, $\overline{\varrho}^{(n)}(t)$, and recall $(X^{(n)}_1(t),\overline{X}^{(n)}_1(t))\stackrel{d}{=}(X^{(n)}_2(t),\overline{X}^{(n)}_2(t))\stackrel{d}{=}\cdots\stackrel{d}{=}(X^{(n)}_n(t),\overline{X}^{(n)}_n(t))$ to control the first probability on the right-hand side of \eqref{2probab_est} by
\begin{equation} \label{eqn: recast1}
\frac{2\sqrt{n}}{\eps}\,\E\Big[\Big|F_{0,M}\big(X^{(n)}_1(t)\big)-F_{0,M}\big(\overline{X}^{(n)}_1(t)\big)\Big|\Big].
\end{equation}
In view of the assumptions on $f_0$ and $g_M$, $M>0$, we can find a constant $C<\infty$ independent of $M$ such that $|f_0(x)|\le Ce^{C|x|}$, $x\in\rr$ and $|g_M(x)|\le C\,\mathbf{1}_{\{|x|>M\}}$, $x\in\rr$, $M>0$. This and the convexity of the absolute value function show that the expression in \eqref{eqn: recast1} is not greater than
\begin{equation}
\begin{split}
& \;\frac{2\sqrt{n}C^2}{\eps}\,\E\Big[e^{C|X^{(n)}_1(t)|}\,\mathbf{1}_{\{|X^{(n)}_1(t)|>M\}}\,\mathbf{1}_{\{|X^{(n)}_1(t)|\ge|\overline{X}^{(n)}_1(t)|\}}\,\big|X^{(n)}_1(t)-\overline{X}^{(n)}_1(t)\big|\Big] \\
& +\frac{2\sqrt{n}C^2}{\eps}\,\E\Big[e^{C|\overline{X}^{(n)}_1(t)|}\,\mathbf{1}_{\{|\overline{X}^{(n)}_1(t)|>M\}}\,\mathbf{1}_{\{|\overline{X}^{(n)}_1(t)|>|X^{(n)}_1(t)|\}}\,\big|X^{(n)}_1(t)-\overline{X}^{(n)}_1(t)\big|\Big]. 
\end{split}
\end{equation}
Leaving out the second indicator random variables from both expectations and applying H\"older's inequality twice we end up with 
\begin{equation}
\begin{split}
& \;\frac{2\sqrt{n}C^2}{\eps}\,\E\Big[e^{3C|X^{(n)}_1(t)|}\Big]^{1/3}\,
\pp\big(\big|X^{(n)}_1(t)\big|>M\big)^{1/3}\,
\E\Big[\big|X^{(n)}_1(t)-\overline{X}^{(n)}_1(t)\big|^3\Big]^{1/3} \\
& +\frac{2\sqrt{n}C^2}{\eps}\,\E\Big[e^{3C|\overline{X}^{(n)}_1(t)|}\Big]^{1/3}\,
\pp\big(\big|\overline{X}^{(n)}_1(t)\big|>M\big)^{1/3}\,
\E\Big[\big|X^{(n)}_1(t)-\overline{X}^{(n)}_1(t)\big|^3\Big]^{1/3},
\end{split}
\end{equation}
which tends to $0$ when one takes the limits superior $n\to\infty$, $M\to\infty$ due to the estimates \eqref{Ecomp}, \eqref{Pcomp} and the $p=3$ version of the inequality \eqref{prop_of_chaos1}.

\medskip

An appeal to Markov's inequality and the independence of $\overline{X}^{(n)}_1(t)\stackrel{d}{=}\overline{X}^{(n)}_2(t)\stackrel{d}{=}\cdots\stackrel{d}{=}\overline{X}^{(n)}_n(t)\stackrel{d}{=}\varrho(t)$ reveal that the second probability on the right-hand side of \eqref{2probab_est} is at most 
\begin{equation}
\frac{4}{\eps^2}\,\E\bigg[\Big(F_{0,M}\big(\overline{X}^{(n)}_1(t)\big)-
\big\langle F_{0,M},\varrho(t)\big\rangle\Big)^{\!2}\bigg].
\end{equation}
Moreover, by the definition of $F_{0,M}$, $M>0$ and the assumptions on $f_0$ and $g_M$, $M>0$ we have $F_{0,M}(x)=0$, $x\in[-M,M]$, $M>0$ and $|F_{0,M}(x)|\le Ce^{C|x|}$, $|x|>M$, $M>0$, which allows to upper bound the latter expectation by
\begin{equation}
C^2\,\E\Big[\mathbf{1}_{\{|\overline{X}^{(n)}_1(t)|>M\}}\,
e^{2C|\overline{X}^{(n)}_1(t)|}\Big]
\le C^2\,\pp\big(\big|\overline{X}^{(n)}_1(t)\big|>M\big)^{1/2}\,
\E\Big[e^{4C|\overline{X}^{(n)}_1(t)|}\Big]^{1/2}.
\end{equation}  
To finish the proof of the lemma we pass to the limits superior $n\to\infty$, $M\to\infty$ relying on the estimates \eqref{Pcomp}, \eqref{Ecomp} one more time. \ep

\medskip

\noindent\textbf{Convergence of $I_2^{(n)}(\cdot)$.} We claim that, for all $t\ge0$, it holds 
\begin{equation} \label{eqn: convI1}
I_2^{(n)}(t)\stackrel{n\to\infty}\longrightarrow
\Big(\big\langle f_1,\varrho(t)\big\rangle,\,\ldots,\,\big\langle f_k,\varrho(t)\big\rangle\Big)
\end{equation}
in probability, which together with \eqref{thm1pf_equ1} and \eqref{thm1pf_equ15} yields Theorem \ref{thm1}. To obtain \eqref{eqn: convI1}, we need to establish $\lim_{n\to\infty} \langle f_j,\widetilde{\varrho}^{(n)}(t)\rangle=\langle f_j,\varrho(t)\rangle$ in probability for every fixed $t\ge0$ and $j\in\{1,2,\ldots,k\}$. Consider the decomposition 
\begin{equation}\label{I1decomp}
\big\langle f_j,\widetilde{\varrho}^{(n)}(t)\big\rangle
=\big\langle f_j,\varrho(t)\big\rangle
+\xi^{(n)}(t)\,\big\langle f_j\,h_M,\varrho^{(n)}(t)-\varrho(t)\big\rangle
+\xi^{(n)}(t)\,\big\langle f_j\,\widehat{h}_M,\varrho^{(n)}(t)-\varrho(t)\big\rangle,
\end{equation}
valid for any $M>0$. We have $\lim_{n\to\infty}\, \xi^{(n)}(t)\,\langle f_j\,h_M,\varrho^{(n)}(t)-\varrho(t)\rangle=0$ in probability due to $|\xi^{(n)}(t)|\le 1$ and Proposition \ref{prop1} (note that $f_j\,h_M$ is continuous and bounded). Finally, $|\xi^{(n)}(t)|\le 1$, integration by parts (in which the boundary terms vanish thanks to the estimate \eqref{Ecomp} and Markov's inequality), the union bound, and Lemma \ref{lemma:unifromIntegrability1} give
\begin{equation}\label{thm1pf_equ10}
\begin{split}
& \;\limsup_{M\to\infty}\,\limsup_{n\to\infty}\;
\pp\Big(\big|\xi^{(n)}(t)\,\big\langle f_j\,\widehat{h}_M,\varrho^{(n)}(t)-\varrho(t)\big\rangle\big|>\eps\Big) \\
& \le \limsup_{M\to\infty}\,\limsup_{n\to\infty}\;\pp\bigg(\bigg|\int_\rr f_j'(x)\,\widehat{h}_M(x)\,\big(F_{\varrho^{(n)}(t)}(x)-R(t,x)\big)\,\mathrm{d}x\bigg|>\eps/2\bigg) \\
&\quad\; +\limsup_{M\to\infty}\,\limsup_{n\to\infty}\;\pp\bigg(\bigg|\int_\rr f_j(x)\,\widehat{h}'_M(x)\,\big(F_{\varrho^{(n)}(t)}(x)-R(t,x)\big)\,\mathrm{d}x\bigg|>\eps/2\bigg)=0
\end{split}
\end{equation}
for all $\eps>0$, so that $\lim_{n\to\infty} \langle f_j,\widetilde{\varrho}^{(n)}(t)\rangle=\langle f_j,\varrho(t)\rangle$ in probability as desired. \ep

\section{Proof of Theorem \ref{thm2}} \label{sec:proofCLT}

It is convenient to introduce the truncated versions $\widehat{\tau}^{(n)}:=\tau^{(n)}\wedge(\tau+1)$, $n\in\nn$ of the hitting times $\tau^{(n)}$, $n\in\nn$. The convergence in distribution of $\sqrt{n}\,(\tau^{(n)}-\tau)$ to a limit is then equivalent to the convergence in distribution of $\sqrt{n}\,(\widehat{\tau}^{(n)}-\tau)$ to the same limit thanks to the following proposition, which is proved further below in this section.

\begin{proposition}\label{thm2pf_prop1}
In the setting of Theorem \ref{thm2}, $\tau^{(n)}\stackrel{n\to\infty}{\longrightarrow}\tau$ in probability. 
\end{proposition}

With the simplified notations
\begin{equation}
Z^{(n)}(\cdot):={\mathcal J}_{J;f_1,\ldots,f_k}\big(\varrho^{(n)}(\cdot)\big)\quad\text{and}\quad Z(\cdot):={\mathcal J}_{J;f_1,\ldots,f_k}\big(\varrho(\cdot)\big),
\end{equation}
our starting point for the proof of the convergence of $\sqrt{n}\,(\widehat{\tau}^{(n)}-\tau)$ is the identity
\begin{equation}\label{thm2_basic_identity}
\mathbf{1}_{\{\tau^{(n)}\le\tau+1\}}\,Z^{(n)}(\widehat{\tau}^{(n)})=\mathbf{1}_{\{\tau^{(n)}\le\tau+1\}}\,Z(\tau).
\end{equation}
The latter stems from the continuity of $Z^{(n)}(\cdot)$ and $Z(\cdot)$: for $Z^{(n)}(\cdot)$, it is a direct consequence of the definitions and, for $Z(\cdot)$, one can write $\langle f_1,\varrho(\cdot)\rangle,\,\ldots,\,\langle f_k,\varrho(\cdot)\rangle$ as $\E[f_1(\overline{X}^{(1)}_1(\cdot))],\,\ldots,\,\E[f_k(\overline{X}^{(1)}_1(\cdot))]$ and conclude by taking the expectation in It\^o's formula and using Fubini's theorem (recall $f_1,\ldots,f_k\in\mathcal{E}_3\subset\mathcal{E}_2$ and the estimate \eqref{Ecomp}). We observe in passing that, for the same reasons in conjunction with the dominated convergence theorem, $Z(\cdot)$ is actually continuously differentiable. 

\medskip

Next, we expand \eqref{thm2_basic_identity} into 
\begin{equation}\label{thm2_expansion}
\begin{split}
& \;\mathbf{1}_{\{\tau^{(n)}\le\tau+1\}}\,\sqrt{n}\,\big(Z(\widehat{\tau}^{(n)})-Z(\tau)\big) \\
& = -\mathbf{1}_{\{\tau^{(n)}\le\tau+1\}}\,\sqrt{n}\,\big(Z^{(n)}(\tau)-Z(\tau)\big) \\
& \quad\, + \mathbf{1}_{\{\tau^{(n)}\le\tau+1\}}\,\sqrt{n}\,\big(Z^{(n)}(\tau)-Z(\tau)\big)  
- \mathbf{1}_{\{\tau^{(n)}\le\tau+1\}}\,\sqrt{n}\,\big(Z^{(n)}(\widehat{\tau}^{(n)})-Z(\widehat{\tau}^{(n)})\big).
\end{split}
\end{equation}
In view of the continuous differentiability of $Z(\cdot)$, the mean value theorem and Proposition \ref{thm2pf_prop1}, the left-hand side of \eqref{thm2_expansion} converges in distribution as $n\to\infty$ if and only if $\sqrt{n}\,(\widehat{\tau}^{(n)}-\tau)$ converges in distribution as $n\to\infty$, and the two limits differ by a factor of $Z'(\tau)\neq 0$ (cf. \eqref{asmp_on_a}). Concurrently, the first line on the right-hand side of \eqref{thm2_expansion} tends to $\sum_{j=1}^k 
{\mathcal J}_{J_{x_j};f_1,\ldots,f_k}(\varrho(\tau))\int_\rr f'_j(x)\,G(\tau,x)\,\mathrm{d}x$ in distribution as $n\to\infty$ by Proposition \ref{thm2pf_prop1} and Theorem \ref{thm1}. 

\medskip

To obtain Theorem \ref{thm2} it now suffices to verify that the second line on the right-hand side of \eqref{thm2_expansion} converges to $0$ in probability as $n\to\infty$. As a result of \eqref{thm1pf_equ1}-\eqref{I2def}, the desired convergence follows from
\begin{equation}
\begin{split}
& \;\sqrt{n}\,\big(Z^{(n)}(\tau)-Z(\tau)\big)
-\sqrt{n}\,\big(Z^{(n)}(\widehat{\tau}^{(n)})-Z(\widehat{\tau}^{(n)})\big)  \\
& = I^{(n)}_1(\tau)\,\nabla J\big(I^{(n)}_2(\tau)\big)
-I^{(n)}_1(\widehat{\tau}^{(n)})
\,\nabla J\big(I^{(n)}_2(\widehat{\tau}^{(n)})\big) \\
& =I^{(n)}_1(\tau)\Big(\nabla J\big(I^{(n)}_2(\tau)\big)
-\nabla J\big(I^{(n)}_2(\widehat{\tau}^{(n)})\big)\Big)
+\big(I^{(n)}_1(\tau)-I^{(n)}_1(\widehat{\tau}^{(n)})\big)\,\nabla J\big(I^{(n)}_2(\widehat{\tau}^{(n)})\big),
\end{split}
\end{equation}
\eqref{thm1pf_equ15}, \eqref{eqn: convI1}, and the next two lemmas. 

\begin{lemma}\label{thm2lemma1}
In the setting of Theorem \ref{thm2}, 
\begin{equation}
\nabla J\big(I^{(n)}_2(\widehat{\tau}^{(n)})\big)
\stackrel{n\to\infty}{\longrightarrow}
\nabla J\Big(\big\langle f_1,\varrho(\tau)\big\rangle,\,\ldots,\,\big\langle f_k,\varrho(\tau)\big\rangle\Big)
\end{equation}
in probability. 
\end{lemma}

\begin{lemma}\label{thm2lemma2}
In the setting of Theorem \ref{thm2}, $I^{(n)}_1(\tau)-I^{(n)}_1(\widehat{\tau}^{(n)})\stackrel{n\to\infty}{\longrightarrow}0$ in probability. 
\end{lemma}

We complete the proof of Theorem \ref{thm2} by establishing Proposition \ref{thm2pf_prop1}, Lemma \ref{thm2lemma1}, and Lemma \ref{thm2lemma2}. 

\medskip

\noindent\textbf{Proof of Proposition \ref{thm2pf_prop1}.} Our proof of the proposition relies on the following lemma that extends the convergence result of Proposition \ref{prop1} to test functions in ${\mathcal E}_0$.

\begin{lemma}\label{thm2pf_lemma1}
Let Assumption \ref{main_asmp} be satisfied. Then, for all $f_0\in {\mathcal E}_0$, $T\ge0$ and $\varepsilon>0$,
\begin{equation}\label{thm2pf_equ1}
\lim_{n\to\infty}\;\pp\Big(\sup_{t\in[0,T]}\,\big|\big\langle f_0,\varrho^{(n)}(t)\big\rangle - \big \langle f_0, \varrho(t) \big \rangle \big| >\varepsilon \Big)=0.
\end{equation} 
\end{lemma}

Given Lemma \ref{thm2pf_lemma1}, the uniform continuity of the function $J$ on compact neighborhoods of the set $\{(\langle f_1,\varrho(t) \rangle ,\,\ldots,\,\langle f_k,\varrho(t)\rangle):\,t\in[0,\tau+1]\}\subset\rr^k$ implies
\begin{equation}\label{thm2pf_equ9}
\lim_{n\to\infty}\;\pp\Big(\sup_{t\in[0,\tau+1]} \big|Z^{(n)}(t)-Z(t)\big|>\varepsilon\Big)=0,\quad \varepsilon>0.
\end{equation} 
Since $Z'(\tau)\neq0$ (cf. \eqref{asmp_on_a}) and $Z'(\cdot)$ is continuous, there exist $[0,1]\ni\upsilon_m\downarrow 0$ such that $Z(\tau+\upsilon_m)>Z(\tau)$ for all $m\in\nn$ if $Z'(\tau)>0$, or $Z(\tau+\upsilon_m)<Z(\tau)$ for all $m\in\nn$ if $Z'(\tau)<0$. Applying \eqref{thm2pf_equ9} with $\varepsilon:=|Z(\tau+\upsilon_m)-Z(\tau)|/2=|Z(\tau+\upsilon_m)-a|/2$ consecutively, we find that 
\begin{equation}
\lim_{n\to\infty}\,\pp\big(\tau^{(n)}>\tau+\upsilon_m\big)=0
\end{equation}
for all $m$. At the same time, for all $\upsilon>0$, we have by \eqref{thm2pf_equ9}:
\begin{equation}
\lim_{n\to\infty}\,\pp\big(\tau^{(n)}\le\tau-\upsilon\big)
\le \lim_{n\to\infty}\,\pp\Big(\sup_{t\in[0,\tau-\upsilon]} \big|Z^{(n)}(t)-Z(t)\big|\ge\min_{t\in[0,\tau-\upsilon]} |a-Z(t)|\Big)=0.
\end{equation}

\smallskip

We conclude the proof of the proposition by showing Lemma \ref{thm2pf_lemma1}. Recalling the auxiliary functions $h_M$, $M>0$ and $\widehat{h}_M$, $M>0$ from the beginning of Section \ref{sec:PfLemma3}, we know from Proposition \ref{prop1} that, for any $M>0$ and $\varepsilon>0$,  
\begin{equation}
\lim_{n\to\infty}\;\pp\Big(\sup_{t\in[0,T]}\,\big|\big\langle f_0\,h_M,\varrho^{(n)}(t)\big\rangle - \big \langle f_0\,h_M, \varrho(t) \big \rangle \big| >\varepsilon \Big)=0.
\end{equation} 
Therefore, it is enough to check that
\begin{equation} \label{thm2pf_equ3}
\limsup_{M\to\infty} \limsup_{n\to\infty}\,\pp\Big(\sup_{t\in[0,T]}\,\big| \big \langle f_0\,\widehat{h}_M,\varrho^{(n)}(t)\big\rangle\big|>\varepsilon \Big)=0, \;\;
\limsup_{M\to\infty} \sup_{t\in[0,T]}\,\big|\big\langle f_0\,\widehat{h}_M,\varrho(t)\big\rangle\big|=0.
\end{equation}

For the first assertion in \eqref{thm2pf_equ3}, we use the definition of $\varrho^{(n)}(t)$, the observation $X^{(n)}_1(\cdot)\stackrel{d}{=}X^{(n)}_2(\cdot)\stackrel{d}{=}\cdots\stackrel{d}{=}X^{(n)}_n(\cdot)$, the estimate $|f_0(x)\,\widehat{h}_M(x)|\le Ce^{C|x|}\,\mathbf{1}_{\{|x|>M\}}$, $x\in\rr$, and the Cauchy-Schwarz inequality to deduce that
\begin{equation}\label{thm2pf_equ5}
\begin{split}
\E\Big[\sup_{t\in[0,T]}\,\big|\big\langle f_0\,\widehat{h}_M,\varrho^{(n)}(t)\big\rangle\big|\Big] 
& = \E\bigg[\sup_{t\in[0,T]}\,\bigg|\frac{1}{n}\sum_{i=1}^n 
f_0\big(X_i^{(n)}(t)\big)\,\widehat{h}_M\big(X_i^{(n)}(t)\big)\bigg|\bigg] \\ 
& \le \E\Big[\sup_{t\in[0,T]}\,\Big|f_0\big(X_1^{(n)}(t)\big)\,\widehat{h}_M\big(X_1^{(n)}(t)\big)\Big|\Big] \\
& \le C\,\E\Big[\sup_{t\in[0,T]}\,\Big(e^{C|X_1^{(n)}(t)|}\,\mathbf{1}_{\{|X_1^{(n)}(t)|>M\}}\Big)\Big] \\
& \le C\,\E\Big[\sup_{t\in[0,T]}\,e^{2C|X_1^{(n)}(t)|}\Big]^{1/2}\,
\E\Big[\sup_{t\in[0,T]}\,\mathbf{1}_{\{|X_1^{(n)}(t)|>M\}}\Big]^{1/2}.
\end{split}
\end{equation}
Since $\sup_{t\in[0,T]} \mathbf{1}_{\{|X_1^{(n)}(t)|>M\}}=\mathbf{1}_{\{\sup_{t\in[0,T]} |X_1^{(n)}(t)|>M\}}$, the first assertion in \eqref{thm2pf_equ3} now follows from Markov's inequality and the estimates \eqref{Ecomp}, \eqref{Pcomp}. 

\smallskip

For the second assertion in \eqref{thm2pf_equ3}, we recall that $\overline{X}_1^{(n)}(t)\stackrel{d}{=}\varrho(t)$, $t\ge0$, allowing us to bound $\sup_{t\in[0,T]} |\langle f_0\,\widehat{h}_M,\varrho(t)\rangle|$ by
\begin{equation}
\begin{split}
& \;\E\Big[\sup_{t\in[0,T]}\;\Big|f_0\big(\overline{X}_1^{(n)}(t)\big)\,\widehat{h}_M\big(\overline{X}_1^{(n)}(t)\big)\Big|\Big] \\
& \le C\,\E\Big[\sup_{t\in[0,T]}\,e^{2C|\overline{X}_1^{(n)}(t)|}\Big]^{1/2}\,\pp\Big(\sup_{t\in[0,T]}\,|\overline{X}_1^{(n)}(t)|>M\Big)^{1/2}
\end{split}
\end{equation}
via the procedure in the last paragraph. The estimates \eqref{Ecomp}, \eqref{Pcomp} yield the result. \ep

\medskip

\noindent\textbf{Proof of Lemma \ref{thm2lemma1}.} In view of the continuity of $\nabla J$, it suffices to show that 
\begin{equation}\label{thm2pf_equ25}
\big\langle f_j,\widetilde{\varrho}^{(n)}(\widehat{\tau}^{(n)})\big\rangle\stackrel{n\to\infty}{\longrightarrow}\langle f_j,\varrho(\tau)\rangle,\quad j=1,\,2,\,\ldots,\,k
\end{equation}
in probability. Since $\widetilde{\varrho}^{(n)}(\widehat{\tau}^{(n)})$ is a convex combination of $\varrho^{(n)}(\widehat{\tau}^{(n)})$ and $\varrho(\widehat{\tau}^{(n)})$, we may swap $\widetilde{\varrho}^{(n)}(\widehat{\tau}^{(n)})$ for $\varrho(\widehat{\tau}^{(n)})$ on the left-hand side of \eqref{thm2pf_equ25} by Lemma \ref{thm2pf_lemma1} with $T:=\tau+1$. Then, Proposition \ref{thm2pf_prop1} and the continuity of $\langle f_j,\varrho(\cdot)\rangle$ give the lemma. \ep

\medskip

\noindent\textbf{Proof of Lemma \ref{thm2lemma2}.} We need to verify that
\begin{equation}\label{thm2lemma2goal}
\sqrt{n}\,\big\langle f_j,\varrho^{(n)}(\tau)-\varrho^{(n)}(\widehat{\tau}^{(n)})\big\rangle-\sqrt{n}\,\big\langle f_j,\varrho(\tau)-\varrho(\widehat{\tau}^{(n)})\big\rangle\stackrel{n\to\infty}{\longrightarrow}0,\quad j=1,\,2,\,\ldots,\,k
\end{equation}
in probability. For a fixed $j\in\{1,2,\ldots,k\}$, we start by establishing the corresponding convergence under the assumption that $f_j\in C^3_c(\rr)\subset{\mathcal E}_3$. 

\medskip

\noindent\textbf{Step 1: convergence \eqref{thm2lemma2goal} for $f_j\in C^3_c(\rr)$.} Inserting the definition of $\varrho^{(n)}(\cdot)$ and applying It\^o's formula we find for the first term in \eqref{thm2lemma2goal}:
\begin{equation}\label{thm2pf_equ29}
\begin{split}
& \;\sqrt{n}\,\big\langle f_j,\varrho^{(n)}(\tau)-\varrho^{(n)}(\widehat{\tau}^{(n)})\big\rangle
= \frac{1}{\sqrt{n}}\sum_{i=1}^n \Big(f_j\big(X^{(n)}_i(\tau)\big)
-f_j\big(X^{(n)}_i(\widehat{\tau}^{(n)})\big)\Big) \\
& = \frac{1}{\sqrt{n}}\,\sum_{i=1}^n \int_{\widehat{\tau}^{(n)}}^\tau 
b\big(F_{\varrho^{(n)}(t)}\big(X_i^{(n)}(t)\big)\big)\,f'_j\big(X^{(n)}_i(t)\big)+\frac{\sigma(F_{\varrho^{(n)}(t)}(X_i^{(n)}(t)))^2}{2}\,f''_j\big(X^{(n)}_i(t)\big)\,\mathrm{d}t \\ 
& \quad\, +\frac{1}{\sqrt{n}}\,\sum_{i=1}^n \int_{\widehat{\tau}^{(n)}}^\tau \sigma\big(F_{\varrho^{(n)}(t)}\big(X_i^{(n)}(t)\big)\big)\,f'_j\big(X^{(n)}_i(t)\big)\,\mathrm{d}B_i^{(n)}(t).
\end{split}
\end{equation}
To simplify the second line in \eqref{thm2pf_equ29} we introduce the discrete antiderivatives 
\begin{equation}
({\mathcal I}_n b)(r):=\frac{1}{n}\,\sum_{i=1}^n b(i/n)\,\mathbf{1}_{\{r\ge i/n\}},\;\;\bigg({\mathcal I}_n \frac{\sigma^2}{2}\bigg)(r):=\frac{1}{n}\,\sum_{i=1}^n \frac{\sigma(i/n)^2}{2}\,\mathbf{1}_{\{r\ge i/n\}},\quad r\in[0,1].
\end{equation}
Since the order statistics $X^{(n)}_{(1)}(t)\le X^{(n)}_{(2)}(t)\le\cdots\le X^{(n)}_{(n)}(t)$ are almost surely distinct for Lebesgue almost every $t\ge0$ by \cite[theorem on p.~439]{Kr} for the function $x\mapsto\sum_{1\le i_1<i_2\le n} \mathbf{1}_{\{x_{i_1}=x_{i_2}\}}$, we can now use summation by parts, the piecewise constant nature of $({\mathcal I}_n b)(F_{\varrho^{(n)}(t)}(\cdot))$, and the convention $X^{(n)}_{(n+1)}(t)=\infty$ to compute
\begin{equation}
\begin{split}
\frac{1}{n}\,\sum_{i=1}^n b\big(F_{\varrho^{(n)}(t)}\big(X_i^{(n)}(t)\big)\big)\,f'_j\big(X^{(n)}_i(t)\big)
=\frac{1}{n}\,\sum_{i=1}^n b(i/n)\,f'_j\big(X^{(n)}_{(i)}(t)\big) \\
=\sum_{i=1}^n \Big(({\mathcal I}_n b)(i/n)-({\mathcal I}_n b)\big((i-1)/n\big)\Big)\,f'_j\big(X^{(n)}_{(i)}(t)\big) \\
=\sum_{i=1}^n ({\mathcal I}_n b)(i/n)\,\Big(f'_j\big(X^{(n)}_{(i)}(t)\big)-f'_j\big(X^{(n)}_{(i+1)}(t)\big)\Big) \\
=-\int_\rr ({\mathcal I}_n b)\big(F_{\varrho^{(n)}(t)}(x)\big)\,f''_j(x)\,\mathrm{d}x.
\end{split}
\end{equation}
Similarly, we see that
\begin{equation}
\frac{1}{n}\,\sum_{i=1}^n \frac{\sigma\big(F_{\varrho^{(n)}(t)}\big(X_i^{(n)}(t)\big)\big)^2}{2}\,f''_j\big(X^{(n)}_i(t)\big)
=-\int_\rr \bigg({\mathcal I}_n \frac{\sigma^2}{2}\bigg)\big(F_{\varrho^{(n)}(t)}(x)\big)\,f'''_j(x)\,\mathrm{d}x.
\end{equation}
Consequently, we arrive at
\begin{equation}\label{discrete_term}
\begin{split}
& \;\sqrt{n}\,\big\langle f_j,\varrho^{(n)}(\tau)-\varrho^{(n)}(\widehat{\tau}^{(n)})\big\rangle \\
& =-\sqrt{n}\int_{\widehat{\tau}^{(n)}}^\tau \int_\rr ({\mathcal I}_n b)\big(F_{\varrho^{(n)}(t)}(x)\big)\,f''_j(x)+\bigg({\mathcal I}_n \frac{\sigma^2}{2}\bigg)\big(F_{\varrho^{(n)}(t)}(x)\big)\,f'''_j(x)\,\mathrm{d}x\,\mathrm{d}t \\
&\quad\,+N^{(n)}(\tau)-N^{(n)}(\widehat{\tau}^{(n)}),
\end{split}
\end{equation}
where
\begin{equation}
N^{(n)}(t):=\frac{1}{\sqrt{n}}\,\sum_{i=1}^n \int_0^t \sigma\big(F_{\varrho^{(n)}(s)}\big(X_i^{(n)}(s)\big)\big)\,f'_j\big(X^{(n)}_i(s)\big)\,\mathrm{d}B_i^{(n)}(s),\quad t\geq 0.
\end{equation}

\smallskip

On the other hand, integration by parts and the notion of a generalized solution for the PDE \eqref{prop1_equ1} (see \cite[Definition 3]{Gi}) imply that
\begin{equation}\label{thm2pf_equ34}
\sqrt{n}\,\big\langle f_j,\varrho(\tau)-\varrho(\widehat{\tau}^{(n)})\big\rangle
=-\sqrt{n}\int_{\widehat{\tau}^{(n)}}^\tau \int_\rr B\big(R(t,x)\big)\,f_j''(x)+\Sigma\big(R(t,x)\big)\,f_j'''(x)\,\mathrm{d}x\,\mathrm{d}t, 
\end{equation}
which can be combined with \eqref{discrete_term} to 
\begin{equation} \label{eqn: lebDecomp}
\begin{split}
&\; \sqrt{n}\,\big\langle f_j,\varrho^{(n)}(\tau)-\varrho^{(n)}(\widehat{\tau}
^{(n)})\big\rangle-\sqrt{n}\,\big\langle f_j,\varrho(\tau)-\varrho(\widehat{\tau}^{(n)})\big\rangle \\
& =-\sqrt{n}\int_{\widehat{\tau}^{(n)}}^\tau \int_\rr \Big(({\mathcal I}_n b)\big(F_{\varrho^{(n)}(t)}(x)\big)-B\big(R(t,x)\big)\Big)\,f''_j(x) \\
& \qquad\qquad\qquad\quad +\bigg(\bigg({\mathcal I}_n \frac{\sigma^2}{2}\bigg)\big(F_{\varrho^{(n)}(t)}(x)\big)-\Sigma\big(R(t,x)\big)\bigg)\,f'''_j(x)\,\mathrm{d}x\,\mathrm{d}t \\
&\quad\, +N^{(n)}(\tau)-N^{(n)}(\widehat{\tau}^{(n)}).
\end{split}
\end{equation}

\smallskip

To prove that the right-hand side of \eqref{eqn: lebDecomp} converges to $0$ in probability we note that the Lipschitz property of $b$, $\frac{\sigma^2}{2}$ (cf. Assumption 1.1(b)) yields 
\begin{equation}
\lim_{n\to\infty} \sqrt{n}\sup_{r\in[0,1]} |({\mathcal I}_n b)(r)-B(r)|=\lim_{n\to\infty} \sqrt{n}\sup_{r\in[0,1]} \bigg|\bigg({\mathcal I}_n\frac{\sigma^2}{2}\bigg)(r)-\Sigma(r)\bigg|=0.
\end{equation}
Since, in addition, $B$, $\Sigma$ are Lipschitz (cf. Assumption 1.1(b)) and $f''_j$, $f'''_j$ are bounded, it suffices to obtain the limits in probability 
\begin{equation}\label{2convs}
\lim_{n\to\infty} \sqrt{n}\int_{\widehat{\tau}^{(n)}}^\tau \int_\rr \big|F_{\varrho^{(n)}(t)}(x)-R(t,x)\big|\,\mathrm{d}x\,\mathrm{d}t=0\;\;\text{and}\;\; \lim_{n\to\infty} \big(N^{(n)}(\tau)-N^{(n)}(\widehat{\tau}^{(n)})\big)=0.
\end{equation} 

\smallskip

For the first convergence in \eqref{2convs}, we recall the representation of the $W_1$-distance in \eqref{W1repr} and apply the triangle inequality for the latter together with Markov's inequality and Fubini's theorem to find, for all $\varepsilon,\varepsilon'>0$,
\begin{equation}
\begin{split}
&\; \pp\bigg(\bigg|\sqrt{n}\int_{\widehat{\tau}^{(n)}}^\tau W_1(\varrho^{(n)}(t),\varrho(t))\,\mathrm{d}t\bigg|>\varepsilon\bigg) \\
& \le \pp\big(|\widehat{\tau}^{(n)}-\tau|>\varepsilon'\big)
+\pp\bigg(\sqrt{n}\int_{\tau-\varepsilon'}^{\tau+\varepsilon'} W_1\big(\varrho^{(n)}(t),\varrho(t)\big)\,\mathrm{d}t>\varepsilon\bigg) \\
& \le \pp\big(|\widehat{\tau}^{(n)}-\tau|>\varepsilon'\big)
+\frac{\sqrt{n}}{\varepsilon}\,\int_{\tau-\varepsilon'}^{\tau+\varepsilon'} \E\big[W_1\big(\varrho^{(n)}(t),\overline{\varrho}^{(n)}(t)\big)\big]+\E\big[W_1\big(\overline{\varrho}^{(n)}(t),\varrho(t)\big)\big]\,\mathrm{d}t.
\end{split}
\end{equation}
In view of Propositions \ref{thm2pf_prop1}, \ref{prop2.2} and \ref{prop_2.3}, this estimate tends to $0$ for all $\varepsilon>0$ when we take $n\to\infty$ and then $\varepsilon'\downarrow0$. 

\medskip

For the second convergence in \eqref{2convs}, we compute the quadratic variation process
\begin{equation}\label{quad_var}
[N^{(n)}](t)=\int_0^t \frac{1}{n}\sum_{i=1}^n \sigma\big(F_{\varrho^{(n)}(s)}\big(X_i^{(n)}(s)\big)\big)^2\,f'_j\big(X^{(n)}_i(s)\big)^2\,\mathrm{d}s,\quad t\ge0,
\end{equation}
bound the resulting integrand by a constant $C<\infty$, and use the martingale representation theorem (see e.g. \cite[Chapter 3, Theorem 4.6 and Problem 4.7]{KaSh}) to conclude
\begin{equation}
\pp\big(|N(\tau)-N(\widehat{\tau}^{(n)})| > \varepsilon \big) 
\leq  \pp\big(|\widehat{\tau}^{(n)}-\tau|>\varepsilon'\big)
+ \pp\Big(\sup_{t\in[0,C\varepsilon']} |B^{(1)}_1(t)|>\varepsilon\Big).
\end{equation}
Thanks to Proposition \ref{thm2pf_prop1}, it is now enough to send $n\to\infty$ followed by $\varepsilon'\downarrow0$. 

\medskip

\noindent\textbf{Step 2: convergence \eqref{thm2lemma2goal} for general $f_j\in\mathcal{E}_3$.} With the functions $h_M$, $M>0$ and $\widehat{h}_M$, $M>0$ introduced at the beginning of Section \ref{sec:PfLemma3}, we decompose the left-hand side of \eqref{thm2lemma2goal} into
\begin{equation}\label{thm2decomp}
\begin{split}
&\; \sqrt{n}\,\big\langle f_j\,h_M,\varrho^{(n)}(\tau)-\varrho^{(n)}(\widehat{\tau}^{(n)})\big\rangle-\sqrt{n}\,\big\langle f_j\,h_M,\varrho(\tau)-\varrho(\widehat{\tau}^{(n)})\big\rangle \\
& +\sqrt{n}\,\big\langle f_j\,\widehat{h}_M,\varrho^{(n)}(\tau)-\varrho^{(n)}(\widehat{\tau}^{(n)})\big\rangle-\sqrt{n}\,\big\langle f_j\,\widehat{h}_M,\varrho(\tau)-\varrho(\widehat{\tau}^{(n)})\big\rangle.
\end{split}
\end{equation}
The first line in \eqref{thm2decomp} converges to $0$ in probability as $n\to\infty$ by Step 1, so we focus on the second line in \eqref{thm2decomp}. To move from $\varrho^{(n)}(\cdot)$ to $\overline{\varrho}^{(n)}(\cdot)$ therein we will prove that, for all $\varepsilon>0$,
\begin{equation}\label{rhotorhobar}
\limsup_{M\to\infty}\,\limsup_{n\to\infty}\;\pp\Big(
\Big|\sqrt{n}\,\big\langle f_j\,\widehat{h}_M,\varrho^{(n)}(\tau)-\varrho^{(n)}(\widehat{\tau}^{(n)})-\overline{\varrho}^{(n)}(\tau)+\overline{\varrho}^{(n)}(\widehat{\tau}^{(n)})\big\rangle\Big|>\varepsilon\Big)=0.
\end{equation}

\smallskip

We recall that $\widehat{\tau}^{(n)}\in[0,\tau+1]$, insert the definitions of $\varrho^{(n)}$, $\overline{\varrho}^{(n)}$, and exploit Markov's inequality to bound the probability in \eqref{rhotorhobar} by
\begin{equation} \label{eqn: unifromIneq} 
\begin{split} 
&\; \frac{2\sqrt{n}}{\varepsilon}\,\E\bigg[\sup_{t\in[0,\tau+1]} \bigg|\frac{1}{n}\,\sum_{i=1}^n 
\Big(\big(f_j\,\widehat{h}_M\big)\big(X^{(n)}_i(t)\big)
\!-\!\big(f_j\,\widehat{h}_M\big)\big(\overline{X}^{(n)}_i(t)\big)\Big)\bigg|\bigg] \\
& \le\frac{2\sqrt{n}}{\varepsilon}\,\E\Big[\sup_{t\in[0,\tau+1]} \Big|\big(f_j\,\widehat{h}_M\big)\big(X^{(n)}_1(t)\big)-\big(f_j\,\widehat{h}_M\big)\big(\overline{X}^{(n)}_1(t)\big)\Big|\Big],
\end{split}
\end{equation}
where we have used $(X^{(n)}_1(\cdot),\overline{X}^{(n)}_1(\cdot))\stackrel{d}{=}(X^{(n)}_2(\cdot),\overline{X}^{(n)}_2(\cdot))\stackrel{d}{=}\cdots\stackrel{d}{=}(X^{(n)}_n(\cdot),\overline{X}^{(n)}_n(\cdot))$. Due to the mean value theorem for $f_j\,\widehat{h}_M$, the inequality $|(f_j\,\widehat{h}_M)'(x)|\le Ce^{C|x|}\,\mathbf{1}_{\{|x|>M\}}$, $x\in\rr$, the convexity of the absolute value function, and H\"older's inequality the right-hand side of \eqref{eqn: unifromIneq} is less or equal to
\begin{equation}
\begin{split}
&\; \frac{2\sqrt{n}}{\varepsilon}\,\E\Big[\sup_{t\in[0,\tau+1]}\big|X^{(n)}_1(t)-\overline{X}^{(n)}_1(t)\big|^3 \Big]^{1/3} \\
&\cdot\bigg(C\,\E\Big[\sup_{t\in[0,\tau+1]} e^{3C|X^{(n)}_1(t)|}\Big]^{1/3}\,\pp\Big(\sup_{t\in[0,\tau+1]} \big|X^{(n)}_1(t)\big|>M\Big)^{1/3} \\
& \quad\; +C\,\E\Big[\sup_{t\in[0,\tau+1]} e^{3C|\overline{X}^{(n)}_1(t)|}\Big]^{1/3}\,\pp\Big(\sup_{t\in[0,\tau+1]} \big|\overline{X}^{(n)}_1(t)\big|>M\Big)^{1/3}\bigg).
\end{split}
\end{equation}
At this point, \eqref{rhotorhobar} becomes a consequence of the inequality \eqref{prop_of_chaos1} with $p=3$ and the estimates \eqref{Ecomp}, \eqref{Pcomp}.

\medskip

With $\A_t:=b(R(t,\cdot))\,\frac{\mathrm{d}}{\mathrm{d}x}+\frac{\sigma(R(t,\cdot))^2}{2}\,\frac{\mathrm{d}^2}{\mathrm{d}x^2}$, $t\ge0$, we may now replace $\varrho^{(n)}(\cdot)$ by $\overline{\varrho}^{(n)}(\cdot)$ in the second line of \eqref{thm2decomp} and deduce by means of It\^o's formula that
\begin{equation}
\begin{split}
\sqrt{n}\,\big\langle f_j\,\widehat{h}_M,\overline{\varrho}^{(n)}(\tau)-\overline{\varrho}^{(n)}(\widehat{\tau}^{(n)})\big\rangle 
=\frac{1}{\sqrt{n}}\,\sum_{i=1}^n \Big(\!\big(f_j\,\widehat{h}_M\big)\big(\overline{X}^{(n)}_i(\tau)\big)
-\big(f_j\,\widehat{h}_M\big)\big(\overline{X}^{(n)}_i(\widehat{\tau}^{(n)})\big)\!\Big) \\
= \frac{1}{\sqrt{n}}\,\sum_{i=1}^n \int_{\widehat{\tau}^{(n)}}^\tau 
\big(\A_t(f_j\,\widehat{h}_M)\big)\big(\overline{X}_i^{(n)}(t)\big)\,\mathrm{d}t+\overline{N}^{(n)}(\tau)-\overline{N}^{(n)}(\widehat{\tau}^{(n)}),
\end{split}
\end{equation}
where
\begin{equation}
\overline{N}^{(n)}(t):=\frac{1}{\sqrt{n}}\,\sum_{i=1}^n \int_0^t \sigma\big(R(s,\overline{X}_i^{(n)}(s))\big)\,
\big(f_j\widehat{h}_M\big)'\big(\overline{X}_i^{(n)}(s)\big)
\,\mathrm{d}B_i^{(n)}(s),\quad t\ge0.
\end{equation}
Next, we take the expectation in It\^o's formula for $(f_j\widehat{h}_M)(\overline{X}^{(1)}_1(\tau))-(f_j\widehat{h}_M)(\overline{X}^{(1)}_1(t))$, $t\ge0$ relying on $\overline{X}^{(1)}_1(t)\stackrel{d}{=}\varrho(t)$, $t\ge0$ and employ Fubini's theorem (note $f_j\widehat{h}_M\in{\mathcal E}_2$ and the estimate \eqref{Ecomp}), followed by an evaluation at $t=\widehat{\tau}^{(n)}$ to get
\begin{equation}\label{exp_in_Ito}
\big\langle f_j\,\widehat{h}_M,\varrho(\tau)-\varrho(\widehat{\tau}^{(n)})\big\rangle 
=\int_{\widehat{\tau}^{(n)}}^\tau \E\Big[\big(\A_t(f_j\,\widehat{h}_M)\big)\big(\overline{X}_1^{(1)}(t)\big)\Big]\,\mathrm{d}t.
\end{equation}

\smallskip

We proceed using the union bound, Markov's inequality, and $\widehat{\tau}^{(n)}\in[0,\tau+1]$:
\begin{equation}\label{var_est}
\begin{split}
& \;\pp\Big(\Big|\sqrt{n}\,\big\langle f_j\,\widehat{h}_M,\overline{\varrho}^{(n)}(\tau)-\overline{\varrho}^{(n)}(\widehat{\tau}^{(n)})\big\rangle-\sqrt{n}\,\big\langle f_j\,\widehat{h}_M,\varrho(\tau)-\varrho(\widehat{\tau}^{(n)})\big\rangle\Big|>\varepsilon\Big) \\
& \le \frac{2}{\varepsilon}\,\E\bigg[\int_0^{\tau+1}\!\bigg|\frac{1}{\sqrt{n}} \sum_{i=1}^n \Big(\big(\A_t(f_j\,\widehat{h}_M)\big)\big(\overline{X}_i^{(n)}(t)\big) 
-\E\Big[\big(\A_t(f_j\,\widehat{h}_M)\big)\big(\overline{X}_i^{(n)}(t)\big)\Big]\Big)\bigg|\,\mathrm{d}t\bigg] \\
&\quad +\pp\Big(\big|\overline{N}^{(n)}(\tau)-\overline{N}^{(n)}(\widehat{\tau}^{(n)})\big|>\varepsilon/2\Big) \\
& \le\frac{2}{\varepsilon} \int_0^{\tau+1} \mathbb{SD}\Big(\big(\A_t(f_j\,\widehat{h}_M)\big)\big(\overline{X}_1^{(1)}(t)\big)\Big)\,\mathrm{d}t 
+\pp\Big(\big|\overline{N}^{(n)}(\tau)-\overline{N}^{(n)}(\widehat{\tau}^{(n)})\big|>\varepsilon/2\Big),
\end{split}
\end{equation} 
where we have applied Fubini's theorem, Jensen's inequality, and the independence of $\overline{X}^{(n)}_1(t)\stackrel{d}{=}\overline{X}^{(n)}_2(t)\stackrel{d}{=}\cdots\stackrel{d}{=}\overline{X}^{(n)}_n(t)\stackrel{d}{=}\varrho(t)$ and have written $\mathbb{SD}$ for the standard deviation operator. Since the standard deviation of a random variable does not exceed its $L^2$-norm, the boundedness of $b$, $\sigma$ (cf. Assumption \ref{main_asmp}(b)), $f_j\in{\mathcal E}_2$, and the properties of $\widehat{h}_M$ imply an estimate of the form
\begin{equation}
\mathbb{SD}\Big(\big(\A_t(f_j\,\widehat{h}_M)\big)\big(\overline{X}_1^{(1)}(t)\big)\Big)
\le\E\Big[Ce^{C|\overline{X}^{(1)}_1(t)|}\,\mathbf{1}_{\{|\overline{X}^{(1)}_1(t)|>M\}}\Big]^{1/2}.
\end{equation}
Moreover, its right-hand side tends to $0$ as $M\to\infty$ uniformly in $t\in[0,\tau]$ by the Cauchy-Schwarz inequality and the estimates \eqref{Ecomp}, \eqref{Pcomp}.

\medskip

To finish the proof we need to analyze the last probability in \eqref{var_est}. For this purpose, we compute 
\begin{equation}
\big[\overline{N}^{(n)}\big](t)=\frac{1}{n}\sum_{i=1}^n \int_0^t \sigma\big(R(s,\overline{X}_i^{(n)}(s))\big)^2\,\big(f_j\,\widehat{h}_M\big)'\big(\overline{X}^{(n)}_i(s)\big)^2\,\mathrm{d}s,\quad t\ge0.
\end{equation}
Thus, the martingale representation theorem (see e.g. \cite[Chapter 3, Theorem 4.6 and Problem 4.7]{KaSh}) and the union bound give, for all $\varepsilon'>0$,
\begin{equation}\label{thm2pf_equ61}
\begin{split}
&\; \pp\Big(\big|\overline{N}^{(n)}(\tau)-\overline{N}^{(n)}(\widehat{\tau}^{(n)})\big|>\varepsilon/2\Big) \\
& \le \pp\Big(\Big|\big[\overline{N}^{(n)}\big](\tau)-\big[\overline{N}^{(n)}\big](\widehat{\tau}^{(n)})\Big|>\varepsilon'\Big)
+\pp\Big(\sup_{t\in[0,\varepsilon']} \big|B^{(1)}_1(t)\big|>\varepsilon/2\Big).
\end{split}
\end{equation}
Due to $\widehat{\tau}^{(n)}\in[0,\tau+1]$, Markov's inequality, $\overline{X}^{(n)}_1(\cdot)\stackrel{d}{=}\overline{X}^{(n)}_2(\cdot)\stackrel{d}{=}\cdots\stackrel{d}{=}\overline{X}^{(n)}_n(\cdot)\stackrel{d}{=}\overline{X}^{(1)}_1(\cdot)$, Fubini's theorem, the boundedness of $\sigma$ (cf. Assumption \ref{main_asmp}(b)), $f_j\in\mathcal{E}_1$, and the properties of $\widehat{h}_M$ we have an estimate of the type
\begin{equation}
\pp\Big(\Big|\big[\overline{N}^{(n)}\big](\tau)-\big[\overline{N}^{(n)}\big](\widehat{\tau}^{(n)})\Big|>\varepsilon'\Big)
\le\frac{1}{\varepsilon'}\,\int_0^{\tau+1} \E\Big[Ce^{C|\overline{X}^{(1)}_1(t)|}\,\mathbf{1}_{\{|\overline{X}^{(1)}_1(t)|>M\}}\Big]\,\mathrm{d}t.  
\end{equation}
The latter converges to $0$ as $M\to\infty$ thanks to the Cauchy-Schwarz inequality and the estimates \eqref{Ecomp}, \eqref{Pcomp}. It remains to observe that the second probability on the right-hand side of \eqref{thm2pf_equ61} vanishes as $\varepsilon'\downarrow0$. \ep

\section{Applications in stochastic portfolio theory}\label{sec:observable}

\subsection{Dynamics of the market diversity}\label{sec:div_dyn}

Consider a stock market with $n$ companies, as described by the market weight processes $\mu_1(\cdot),\,\mu_2(\cdot),\,\ldots,\,\mu_n(\cdot)$, i.e. the fractions of the total market capital invested in the different companies at any given time. In this context, a concept that has attracted much interest, both for scientific reasons and its importance in investment decisions, is the \textit{market diversity}. Informally speaking, a market is thought of as diverse when one can be certain that no single company will end up with the vast majority of the market capital. In \cite{Fe}, \textsc{Fernholz} has proposed to formalize the notion of diversity as follows. 

\begin{definition}[\cite{Fe}, Definition 2.2.1]\label{def:div}
A market is called \emph{diverse} if for some $\varepsilon>0$ it holds $\max_{1\le i\le n} \mu_i(t)\le 1-\varepsilon$ for all $t\ge0$ almost surely. A market is referred to as \emph{weakly diverse} on a finite time interval $[0,T]$ if for some $\varepsilon>0$ one has 
\begin{equation}
\dfrac{1}{T}\int_0^T \max_{1\le i\le n} \mu_i(t)\,\mathrm{d}t\leq 1-\varepsilon
\end{equation}
almost surely.
\end{definition}

Subsequently, it is noticed in \cite{Fe} that the vector of the market weight processes $\mu(\cdot):=(\mu_1(\cdot),\mu_2(\cdot),\ldots,\mu_n(\cdot))$ takes values in the closed unit simplex 
\begin{equation}
\overline{\Delta}^n:=\bigg\{x\in[0,1]^n:\;\sum_{i=1}^n x_i=1\bigg\},
\end{equation}
whereas the diversity condition $\max_{1\le i\le n} \mu_i(t)\le 1-\varepsilon$ is violated when $\mu(\cdot)$ enters the corresponding open neighborhoods of the vertices of $\overline{\Delta}^n$. Hence, it is natural to use a symmetric concave function on $\overline{\Delta}^n$, which necessarily attains its minimum at the vertices, to quantify the diversity of a market (or the lack thereof). The main examples of such functions discussed in \cite{Fe} are: 
\begin{enumerate}[(i)]
\item the entropy function $H(t)=-\sum_{i=1}^n \mu_i(t)\log\mu_i(t)$, $t\ge0$,
\item the $\ell^p$-norms $D_p(t)=\big(\sum_{i=1}^n \mu_i(t)^p\big)^{1/p}$, $t\ge0$ for $p\in(0,1)$,
\item and the geometric mean $S(t)=\big(\prod_{i=1}^n \mu_i(t)\big)^{1/n}$, $t\ge0$.  
\end{enumerate}
In particular, the entropy function and the $\ell^p$-norms for $p\in(0,1)$ can be employed to test if a market is diverse in the sense of Definition \ref{def:div} (cf. \cite[Proposition 2.3.2]{Fe}). 

\begin{proposition}
A market is diverse if and only if for some $\varepsilon'>0$ it holds $H(t)\ge\varepsilon'$ for all $t\ge0$ almost surely or, equivalently, for some $p\in(0,1)$ and $\varepsilon''>0$ one has $D_p(t)\ge 1+\varepsilon''$ for all $t\ge0$ almost surely.
\end{proposition}

Our Theorem \ref{thm1} can be utilized to capture the dynamics of the entropy $H(\cdot)$, the $\ell^p$-norms $D_p(\cdot)$, $p\in(0,1)$, and the geometric mean $S(\cdot)$ in rank-based models with a large number $n$ of companies. In that setting, the market weight processes are defined in terms of the solution to \eqref{intro_equ1} by 
\begin{equation} \label{eqn:mui}
\mu_i(\cdot)=\frac{e^{X_i^{(n)}(\cdot)}}{e^{X_1^{(n)}(\cdot)}+e^{X_2^{(n)}(\cdot)}+\cdots+e^{X_n^{(n)}(\cdot)}},\quad i=1,\,2,\,\ldots,\,n
\end{equation}
and give rise to the associated entropy, $\ell^p$-norm and geometric mean processes via the items (i), (ii), (iii) above.

\begin{corollary}\label{cor1}
Under Assumption \ref{main_asmp} the following convergences hold in the finite-dimensional distribution sense:
\begin{enumerate}[(a)]
\item for the entropy process $H(\cdot)$,
\begin{equation}
\begin{split}
& \sqrt{n}\,\bigg(H(\cdot)-\log n-\log\big\langle e^x,\varrho(\cdot)\big\rangle+\frac{\langle xe^x,\varrho(\cdot)\rangle}{\langle e^x,\varrho(\cdot)\rangle}\bigg) \\
& \stackrel{n\to\infty}{\longrightarrow}
-\bigg(\frac{1}{\langle e^x,\varrho(\cdot)\rangle}+\frac{\langle xe^x,\varrho(\cdot)\rangle}{\langle e^x,\varrho(\cdot)\rangle^2}\bigg)\int_\rr e^x\,G(\cdot,x)\,\mathrm{d}x-\frac{1}{\langle e^x,\varrho(\cdot)\rangle}\int_\rr \big(e^x+xe^x\big)\,G(\cdot,x)\,\mathrm{d}x,
\end{split}
\end{equation}
\item for an $\ell^p$-norm process $D_p(\cdot)$ with $p\in(0,1)$,
\begin{equation}
\begin{split}
& \sqrt{n}\,\bigg(n^{\frac{p-1}{p}} D_p(\cdot)-\frac{\langle e^{px},\varrho(\cdot)\rangle^{1/p}}{\langle e^x,\varrho(\cdot)\rangle}\bigg) \\
& \stackrel{n\to\infty}{\longrightarrow}
-\frac{\langle e^{px},\varrho(\cdot)\rangle^{1/p-1}}{\langle e^x,\varrho(\cdot)\rangle}\int_\rr e^{px}\,G(\cdot,x)\,\mathrm{d}x
+\frac{\langle e^{px},\varrho(\cdot)\rangle^{1/p}}{\langle e^x,\varrho(\cdot)\rangle^2}\int_\rr e^x\,G(\cdot,x)\,\mathrm{d}x,
\end{split}
\end{equation}
\item for the the geometric mean process $S(\cdot)$,
\begin{equation}
\sqrt{n}\,\bigg(nS(\cdot)-\frac{e^{\langle x,\varrho(\cdot)\rangle}}{\langle e^x,\varrho(\cdot)\rangle}\bigg)
\stackrel{n\to\infty}{\longrightarrow} -\frac{e^{\langle x,\varrho(\cdot)\rangle}}{\langle e^x,\varrho(\cdot)\rangle}\int_\rr G(\cdot,x)\,\mathrm{d}x + \frac{e^{\langle x,\varrho(\cdot)\rangle}}{\langle e^x,\varrho(\cdot)\rangle^2}\int_\rr e^x\,G(\cdot,x)\,\mathrm{d}x.
\end{equation}
\end{enumerate}
\end{corollary}

\noindent\textbf{Proof.} The corollary is a direct consequence of Theorem \ref{thm1}. For the sake of completeness, we write out the functions $J$, $f_1,\,\ldots,\,f_k$ in each of the three cases. 
\begin{enumerate}[(a)]
\item For the normalized entropy process $H(\cdot)-\log n$, take 
\begin{equation}\label{cor1_equ4}
J:\;(0,\infty)\times\rr\to\rr,\;\;(x_1,x_2)\mapsto \log x_1-\frac{x_2}{x_1}, \quad f_1(x)=e^x, \quad f_2(x)=xe^x.
\end{equation}
\item For every normalized $\ell^p$-norm process $n^{\frac{p-1}{p}}D_p(\cdot)$, define
\begin{equation}\label{cor1_equ5}
J:\;(0,\infty)\times(0,\infty)\to\rr,\;\;(x_1,x_2)\mapsto\frac{x_1^{1/p}}{x_2},\quad f_1(x)=e^{px},\quad f_2(x)=e^x.
\end{equation}
\item For the normalized geometric mean process $\frac{S(\cdot)}{n}$, pick 
\begin{equation}\label{cor1_equ6}
J:\;\rr\times(0,\infty)\to\rr,\;\;(x_1,x_2)\mapsto\frac{e^{x_1}}{x_2},\quad
f_1(x)=x,\quad f_2(x)=e^x.
\end{equation}
\end{enumerate}
It is elementary to check the assumptions of Theorem \ref{thm1} for all of these functions. \ep

\subsection{Hitting times of the market diversity}

In this subsection, we discuss the implications of Theorem \ref{thm2} for the measures of diversity from Subsection \ref{sec:div_dyn}. To this end, we denote by $H^*(\cdot)$, $D^*_p(\cdot)$, $p\in(0,1)$, and $S^*(\cdot)$ the limiting entropy, $\ell^p$-norms, and geometric mean processes, respectively: 
\begin{equation}\label{cor2_equ1}
H^*(\cdot)=\log\big\langle e^x,\varrho(\cdot)\big\rangle-\frac{\langle xe^x,\varrho(\cdot)\rangle}{\langle e^x,\varrho(\cdot)\rangle},\;
D_p^*(\cdot)=\frac{\langle e^{px},\varrho(\cdot)\rangle^{1/p}}{\langle e^x,\varrho(\cdot)\rangle},\,p\in(0,1),\;
S^*(t)=\frac{e^{\langle x,\varrho(\cdot)\rangle}}{\langle e^x,\varrho(\cdot)\rangle}.
\end{equation} 
The functions $J$, $f_1,\,\ldots,\,f_k$ in the three cases, given explicitly in \eqref{cor1_equ4}, \eqref{cor1_equ5}, and \eqref{cor1_equ6}, respectively, satisfy the assumptions in Theorem \ref{thm2}, so that one only needs to verify the statements in \eqref{asmp_on_a} for the coefficients $b$, $\sigma$ and levels $a$ of interest. The next proposition provides the dynamics of $H^*(\cdot)$, $D^*_p(\cdot)$, $p\in(0,1)$, and $S^*(\cdot)$, thus, yielding a sufficient condition on $b$, $\sigma$, and $a$ for Theorem \ref{thm2} to apply. 
 
\begin{proposition} \label{prop:driftODE}
Under Assumption \ref{main_asmp} consider 
\begin{equation}
\varrho_p(\cdot):=\frac{e^{px}\,\varrho(\cdot)}{\langle e^{px},\varrho(\cdot)\rangle}\in C([0,\infty),M_1(\rr)), \quad p\in[0,1).
\end{equation}
Then, one has for the processes $H^*(\cdot)$, $D_p^*(\cdot)$, $p\in(0,1)$, and $S^*(\cdot)$ of \eqref{cor2_equ1}:
\begin{eqnarray}
&&\quad\; \frac{\mathrm{d}H^*(t)}{\mathrm{d}t} =-\frac{1}{2}\big\langle \sigma(R(t,\cdot))^2,\varrho_1(t)\big\rangle
-\mathrm{cov}_{\varrho_1(t)}\bigg(x,\,b(R(t,\cdot))+\frac{\sigma(R(t,\cdot))^2}{2}\bigg), \label{cor2_equ2} \\
&&\quad\; \frac{\mathrm{d}D_p^*(t)}{\mathrm{d}t} \!=\! D_p^*(t)\bigg(\!\!\bigg\langle\!\! b(R(t,\cdot)\!)\!+\!\frac{p\,\sigma(R(t,\cdot)\!)^2}{2},\varrho_p(t)\!\!\bigg\rangle
\!\!-\!\bigg\langle\!\! b(R(t,\cdot)\!)\!+\!\frac{\sigma(R(t,\cdot)\!)^2}{2},\varrho_1(t)\!\!\bigg\rangle\!\!\bigg), \label{cor2_equ3} \\
&&\quad\; \frac{\mathrm{d}S^*(t)}{\mathrm{d}t} = S^*(t)\,\bigg(\big\langle b(R(t,\cdot)),\varrho(t)\big\rangle-\bigg\langle b(R(t,\cdot))+\frac{\sigma(R(t,\cdot))^2}{2},\varrho_1(t)\bigg\rangle\bigg). \label{cor2_equ4}
\end{eqnarray}
In particular, whenever $b+\frac{\sigma^2}{2}$ is an increasing function and $a\in(-\infty,H^*(0)]$, $a\in(0,D^*_p(0)]$, $p\in(0,1)$, or $a\in(0, S^*(0)]$, the assertions in \eqref{asmp_on_a} hold and, hence, also the conclusion of Theorem \ref{thm2} for the resulting hitting times of the normalized processes $H(\cdot)-\log n$, $n^{\frac{p-1}{p}}D_p(\cdot)$, $p\in(0,1)$, or $nS(\cdot)$. 
\end{proposition}

We prepare the following continuous version of Chebyshev's sum inequality for the proof of Proposition \ref{prop:driftODE}. 

\begin{lemma}\label{lem5.7}
For all $\nu\in M_1(\rr)$ and increasing functions $f,g$ on $\rr$ integrable with respect to $\nu$,
\begin{equation}\label{eqn:ineqIncLambda}
\langle fg,\nu\rangle\ge\langle f,\nu\rangle\,\langle g,\nu\rangle.
\end{equation}  
\end{lemma}

\noindent\textbf{Proof.} Since $g$ is increasing, there exists an $x_0\in\rr$ such that $g(x)\le\langle g,\nu\rangle$ if $x<x_0$ and $g(x)\ge\langle g,\nu\rangle$ if $x>x_0$. By distinguishing between $x<x_0$ and $x\ge x_0$ and using that $f$ is increasing we deduce
\begin{equation}
f(x)\big(g(x)-\langle g,\nu\rangle\big)\ge f(x_0)\big(g(x)-\langle g,\nu\rangle\big),\quad x\in\rr. 
\end{equation}
Integrating both sides with respect to $\nu$ and rearranging we arrive at \eqref{eqn:ineqIncLambda}. \ep

\medskip

We are now ready to present the proof of Proposition \ref{prop:driftODE}.

\medskip

\noindent\textbf{Proof of Proposition \ref{prop:driftODE}.} We recall the notation $\A_t=b(R(t,\cdot))\,\frac{\mathrm{d}}{\mathrm{d}x}+\frac{\sigma(R(t,\cdot))^2}{2}\,\frac{\mathrm{d}^2}{\mathrm{d}x^2}$, $t\ge0$ and that, for any $f\in{\mathcal E}_2$, 
\begin{equation}
\big\langle f,\varrho(t)\big\rangle-\big\langle f,\varrho(0)\big\rangle
=\int_0^t \big\langle \A_s f,\varrho(s)\big\rangle\,\mathrm{d}s,\quad t\ge0
\end{equation} 
(cf. \eqref{exp_in_Ito}). Due to the boundedness of $b$, $\sigma$ (cf. Assumption \ref{main_asmp}(b)) and $f\in{\mathcal E}_2$ the dominated convergence theorem implies that the function $s\mapsto\langle \A_s f,\varrho(s)\rangle$ is continuous on $[0,\infty)$ and, thus,
\begin{equation}
\frac{\mathrm{d}\langle f,\varrho(t)\rangle}{\mathrm{d}t}=\big\langle\A_t f,\varrho(t)\big\rangle,\quad t\ge0. 
\end{equation}
Therefore, in the setting of Theorem \ref{thm2},
\begin{equation}\label{chain_rule}
\frac{\mathrm{d}\mathcal{J}_{J;f_1,\ldots,f_k}(\varrho(t))}{\mathrm{d}t}=\sum_{j=1}^k {\mathcal J}_{J_{x_j};f_1,\ldots,f_k}\big(\varrho(t)\big)\,\big\langle \A_t f_j,\varrho(t)\big\rangle,\quad t\ge0.
\end{equation}
To obtain the differential equations \eqref{cor2_equ2}, \eqref{cor2_equ3}, and \eqref{cor2_equ4} it suffices to insert into \eqref{chain_rule} the formulas from \eqref{cor1_equ4}, \eqref{cor1_equ5}, and \eqref{cor1_equ6}, respectively, and to simplify the result. 

\medskip

Supposing, in addition, that $b+\frac{\sigma^2}{2}$ is increasing we can first employ Lemma \ref{lem5.7} with $\nu=\varrho_1(t)$, $f(x)=x$ and $g(x)=b(R(t,x))+\frac{\sigma(R(t,x))^2}{2}$ to find 
\begin{equation}
\mathrm{cov}_{\varrho_1(t)}\bigg(x,\,b(R(t,\cdot))+\frac{\sigma(R(t,\cdot))^2}{2}\bigg)\ge0,\quad t\ge0. 
\end{equation}
Consequently, we read off from \eqref{cor2_equ2} that 
\begin{equation}
\frac{\mathrm{d}H^*(t)}{\mathrm{d}t}\le-\frac{1}{2}\,\min_{r\in[0,1]} \sigma(r)^2,\quad t\ge0,
\end{equation}
so \eqref{asmp_on_a} must hold for all $a\in(-\infty,H^*(0)]$, and the conclusion of Theorem \ref{thm2} applies to the hitting times of such $a$ by $H(\cdot)-\log n$.

\medskip

Now, we take $\nu=\varrho_p(t)$, $f(x)=b(R(t,x))+\frac{\sigma(R(t,x))^2}{2}$ and $g(x)=e^{(1-p)x}$, with $p\in[0,1)$, in Lemma \ref{lem5.7} to get
\begin{equation}
\bigg\langle b(R(t,\cdot))+\frac{\sigma(R(t,\cdot))^2}{2},\varrho_1(t)\bigg\rangle\ge 
\bigg\langle b(R(t,\cdot))+\frac{\sigma(R(t,\cdot))^2}{2},\varrho_p(t)\bigg\rangle,\quad 
t\ge0.
\end{equation}
The values of $p\in(0,1)$ and $p=0$ reveal
\begin{equation}
\frac{\mathrm{d}\log D_p^*(\cdot)}{\mathrm{d}t} \leq -\dfrac{1-p}{2}\,\min_{r\in[0,1]} \sigma(r)^2,\;\;p\in(0,1) \quad \text{and} \quad 
\frac{\mathrm{d}\log S^*(\cdot)}{\mathrm{d}t} \leq -\frac{1}{2}\,\min_{r\in[0,1]} \sigma(r)^2,
\end{equation}
respectively, yielding the remaining assertions. \ep

\begin{rmk}
A verification of the conditions in \eqref{asmp_on_a} beyond the setup in Proposition \ref{prop:driftODE} seems to require information on $\varrho(\cdot)$ or, equivalently, $R$ that needs to be deduced on a case-by-case basis. This is possible, for example, when \eqref{intro_equ1} is of the special form
\begin{equation}
\mathrm{d}X_i^{(n)}(t)=\big(2C_1F_{\varrho^{(n)}(t)}\big(X_i^{(n)}(t)\big)+C_2\big)\,\mathrm{d}t+\sigma\,\mathrm{d}B^{(n)}_i(t),\quad i=1,\,2,\,\ldots,\,n
\end{equation}
for some $C_1\neq 0$, $C_2\in\rr$ and $\sigma>0$. Indeed, then the Cauchy problem for the porous medium equation \eqref{prop1_equ1} reduces to the one for the generalized Burgers equation
\begin{equation}
R_t=-(C_1 R^2+C_2R)_x+\frac{\sigma^2}{2}\,R_{xx},\quad R(0,\cdot)=F_{\lambda}(\cdot).
\end{equation}
The solution of the latter is provided by the Cole-Hopf transformation $R=-\frac{2C_1}{\sigma^2}(\log\varphi)_x$, where $\varphi$ is the solution of the Cauchy problem for the heat equation
\begin{equation}
\varphi_t=-C_2\,\varphi_x+\frac{\sigma^2}{2}\,\varphi_{xx},\quad \varphi(0,x)=e^{-\frac{\sigma^2}{2C_1}\int_0^x F_\lambda(y)\,\mathrm{d}y}.
\end{equation}
For any fixed $\lambda\in M_1(\rr)$ (perhaps retrieved from the observed market capitalizations), $\varphi$ is given explicitly by a convolution with the heat kernel, and one can check if the conditions in \eqref{asmp_on_a} are valid for the resulting $\varrho(\cdot)=R_x(\cdot,x)\,\mathrm{d}x$. 
\end{rmk}

\subsection{Performance of functionally generated portfolios} \label{sec:functionallyGenerated}

This last subsection is devoted to a discussion of the performance of multiplicatively and additively generated portfolios $\pi^{\widetilde{\Psi};\times}$ and $\pi^{\widetilde{\Psi};+}$, as defined in the introduction. We focus initially on their associated non-decreasing \emph{excess growth} processes 
\begin{equation}
-\frac{1}{2}\sum_{i,j=1}^n \int_0^t \frac{\widetilde{\Psi}_{x_ix_j}(\mu(\cdot))}{\widetilde{\Psi}(\mu(\cdot))}\,\mathrm{d}[\mu_i,\mu_j](\cdot),\,t\ge0\;\text{and}\;-\frac{1}{2}\sum_{i,j=1}^n \int_0^t \widetilde{\Psi}_{x_ix_j}\big(\mu(\cdot)\big)\,\mathrm{d}[\mu_i,\mu_j](\cdot),\,t\ge0
\end{equation}
which enter the value processes $V^{\widetilde{\Psi};\times}(\cdot)$ and $V^{\widetilde{\Psi};+}(\cdot)$ relative to that of the market portfolio $\mu(\cdot)$ according to \eqref{master_mult} and \eqref{master_add}, respectively. Under the assumptions in \eqref{extra_asmp}, as well as \eqref{mult_range} or \eqref{add_range}, respectively, the former obey the concentration of measure estimate from \cite[Corollary 8]{IPS} (note the symmetry of $\widetilde{\Psi}$ due to \eqref{eqn:PsitildePsi}, \eqref{eqn:JPsi}; the strong law of large numbers for the process $\mu(\cdot)$ in \cite[equation (4.5)]{BFK}; and that the
derivation of \cite[Corollary 8]{IPS} for multiplicatively generated portfolios carries over mutatis mutandis to the case of additive generation).

\begin{proposition}\label{prop5.7}
Suppose the assumptions in \eqref{extra_asmp}, as well as \eqref{mult_range} or \eqref{add_range} are satisfied for some $n\in\nn$. Then, for all $r,t,\varepsilon>0$ and in the notation of Corollary \ref{main_cor},
\begin{equation}\label{prop5.7_equ1}
\mathbb{P}\bigg(-\frac{1}{2t}\,\sum_{i,j=1}^n \int_0^t \frac{\widetilde{\Psi}_{x_ix_j}(\mu(\cdot))}{\widetilde{\Psi}(\mu(\cdot))}\,\mathrm{d}[\mu_i,\mu_j](\cdot)\le r^\times-r\bigg)
\le \bigg\|\frac{\mathrm{d}\kappa^{(n)}}{\mathrm{d}\zeta^{(n)}}
\bigg\|_{L^2(\zeta^{(n)})} e^{-c^\times(r,\eps) t}
\end{equation}
or
\begin{equation}
\,\mathbb{P}\bigg(-\frac{1}{2t}\,\sum_{i,j=1}^n \int_0^t \widetilde{\Psi}_{x_ix_j}\big(\mu(\cdot)\big)\,\mathrm{d}[\mu_i,\mu_j](\cdot)\le r^+-r\bigg)
\le \bigg\|\frac{\mathrm{d}\kappa^{(n)}}{\mathrm{d}\zeta^{(n)}}
\bigg\|_{L^2(\zeta^{(n)})} e^{-c^+(r,\eps) t},
\end{equation}
respectively.
\end{proposition}

If, in addition, Assumption \ref{main_asmp} holds, one can combine Proposition \ref{prop5.7} with Theorem \ref{thm1} by using the union bound and obtain, for all $r,s,t,\eps>0$, the performance estimates
\begin{equation}\label{perf_est1}
\begin{split}
&\;\pp\bigg(V^{\widetilde{\Psi};\times}(t)\le\frac{\mathcal{J}_{J;f_1,\ldots,f_k}(\varrho(t))-s/\sqrt{n}}{\widetilde{\Psi}(\mu(0))}\,e^{(r^\times-r)t}\bigg)\qquad\qquad\quad\;\; \\
& \le\overline{\Phi}(s/\chi_t)\big(1+o_n(1)\big)+
\bigg\|\frac{\mathrm{d}\kappa^{(n)}}{\mathrm{d}\zeta^{(n)}}
\bigg\|_{L^2(\zeta^{(n)})} e^{-c^\times(r,\eps) t}
\end{split}
\end{equation}
or 
\begin{equation}\label{perf_est2}
\begin{split}
&\;\pp\Big(V^{\widetilde{\Psi};+}(t)\le 1+\mathcal{J}_{J;f_1,\ldots,f_k}\big(\varrho(t)\big)-s/\sqrt{n}-\widetilde{\Psi}\big(\mu(0)\big)+(r^+-r)t\Big) \\
& \le \overline{\Phi}(s/\chi_t)\big(1+o_n(1)\big)+ \bigg\|\frac{\mathrm{d}\kappa^{(n)}}{\mathrm{d}\zeta^{(n)}}
\bigg\|_{L^2(\zeta^{(n)})} e^{-c^+(r,\eps) t},
\end{split}
\end{equation}
respectively, where $\overline{\Phi}$ is the standard normal tail cumulative distribution function, $\chi_t$ is the standard deviation of the time $t$ value of the Gaussian process on the right-hand side of \eqref{thm1_result}, and $o_n(1)$ is a quantity tending to $0$ as $n\to\infty$. Complementary to the performance estimates \eqref{perf_est1}, \eqref{perf_est2} for \textit{fixed} times, Corollary \ref{main_cor}, which is proved next, provides a bound on the \textit{random} time it takes for a multiplicatively or additively generated portfolio to reach the desired performance. 

\medskip

\noindent\textbf{Proof of Corollary \ref{main_cor}.} We only give the proof of \eqref{mult_port_bnd}, as \eqref{add_port_bnd} can be shown in the same way. Our starting point is the observation that
\begin{equation}\label{corpf_equ1}
\pp\big(\eta^{\widetilde{\Psi};\times}\ge\tau+s/\sqrt{n}\big)
\le\pp\big(\eta^{\widetilde{\Psi};\times}\ge \tau^{(n)}\big)+\pp\big(\tau^{(n)}\ge \tau+s/\sqrt{n}\big).
\end{equation}
By the definition of $\eta^{\widetilde{\Psi};\times}$, the first of the latter two summands is less or equal to
\begin{equation}
\begin{split}\label{perf_rand_time}
&\;\pp\bigg(V^{\widetilde{\Psi};\times}(\tau^{(n)})\le\frac{a}{\widetilde{\Psi}(\mu(0))}\,e^{(r^\times-r)(\tau-s/\sqrt{n})}\bigg)
\\
&=\pp\bigg(\!-\frac{1}{2}\sum_{i,j=1}^n \int_0^{\tau^{(n)}} \frac{\widetilde{\Psi}_{x_ix_j}(\mu(\cdot))}{\widetilde{\Psi}(\mu(\cdot))}\,\mathrm{d}[\mu_i,\mu_j](\cdot)\le(r^\times-r)\big(\tau-s/\sqrt{n}\big)\!\bigg)
\end{split}
\end{equation}
(recall \eqref{master_mult} and $\widetilde{\Psi}(\mu(\tau^{(n)}))=a$). Since the excess growth process is non-decreasing, the probability on the right-hand side of \eqref{perf_rand_time} is at most
\begin{equation}\label{eq5.31}
\pp\big(\tau^{(n)}\le \tau-s/\sqrt{n}\big)+\pp\bigg(\!-\frac{1}{2}\sum_{i,j=1}^n \int_0^{\tau-s/\sqrt{n}} \frac{\widetilde{\Psi}_{x_ix_j}(\mu(\cdot))}{\widetilde{\Psi}(\mu(\cdot))}\,\mathrm{d}[\mu_i,\mu_j](\cdot)\le(r^\times-r)\big(\tau-s/\sqrt{n}\big)\!\bigg).
\end{equation}
Using \eqref{thm2_result} for the second summand on the right-hand side of \eqref{corpf_equ1} and the first summand in \eqref{eq5.31}, then \eqref{prop5.7_equ1} for the second summand in \eqref{eq5.31} we get \eqref{mult_port_bnd}. \ep

\bibliographystyle{plain}
\bibliography{StocPortfoliloCLT}

\bigskip\medskip

\end{document}